\documentclass[11pt]{article}
\usepackage{amsmath}
\usepackage{amssymb}

\usepackage{multicol}

\usepackage{color}

\newcommand{\la}{\lambda}

{\catcode `\@=11 \global\let\AddToReset=\@addtoreset}
\AddToReset{equation}{section}

\usepackage{epsfig}
\usepackage{amsmath}
\usepackage{amssymb}
\usepackage{graphicx}
\usepackage{color}

\usepackage[latin1]{inputenc}

\newtheorem{proposition}{Proposition}[section]
\newtheorem{@definition}{\sc Definition}[section]

\newtheorem{@remark}{\sc Remark}[section]
\newenvironment{remark}{\begin{@remark}\rm}{\end{@remark}}
\newtheorem{@example}{\sc Example}[section]

\newcommand{\beqn}{\begin{displaymath}}
\newcommand{\eeqn}{\end{displaymath}}
\newcommand{\beq}{\begin{equation}}  
\newcommand{\eeq}{\end{equation}}

\def\mathsf{\bf}

\def\R{\mathbb{R}}
\def\Z{\mathbb{Z}}

\def\i{\mathrm i}
\def\d{\mathrm d}
\def\e{\mathrm e}

\def\E{\mathrm E}
\def\P{\mathrm P}

\def\text{\mbox}

\def\1{{\bf 1}}

\newcommand{\mbf}[1]{\mbox{\boldmath $#1$}}

\newcommand{\nn}{\nonumber}
\newcommand{\noi}{\noindent}

\newcommand{\cY}{\mathcal Y}

\newcommand{\mbt}{\boldsymbol{t}}

\newcommand{\mbx}{\boldsymbol{x}}

\newcommand{\mbgamma}{\boldsymbol{\gamma}}

\newcommand{\smbgamma} {\boldsymbol{\gamma}}

\def\eq2{
\stackrel{\small \rm mod \,2}{=}}

\def\n2{
\stackrel{\small \rm mod \,2}{\neq}}

\def\limd{\renewcommand{\arraystretch}{0.5}
\begin{array}[t]{c}
\stackrel{\rm d}{\longrightarrow} \\
\end{array}\renewcommand{\arraystretch}{1}}

\def\limfdd{\renewcommand{\arraystretch}{0.5}
\begin{array}[t]{c}
\stackrel{\rm fdd}{\longrightarrow} \\
\end{array}\renewcommand{\arraystretch}{1}}

\def\eqfdd{\renewcommand{\arraystretch}{0.5}
\begin{array}[t]{c}
\stackrel{\rm fdd}{=} \\
\end{array}\renewcommand{\arraystretch}{1}}

\def\neqfdd{\renewcommand{\arraystretch}{0.5}
\begin{array}[t]{c}
\stackrel{\rm fdd}{\neq} \\
\end{array}\renewcommand{\arraystretch}{1}}

\newtheorem{thm}{Theorem}[section]

\newtheorem{lem}[thm]{Lemma}
\newtheorem{prop}[thm]{Proposition}

\def\vep{\varepsilon}

\begin{document}

\title{Scaling transition and edge effects for negatively
dependent \\ linear random fields on ${\mathbb{Z}}^2$}

\author{
Donatas Surgailis    \footnote{E-mail: donatas.surgailis@mif.vu.lt}\\
\small \it Vilnius University, Faculty of Mathematics and Informatics, Naugarduko 24, 03225 Vilnius, Lithuania
}
\maketitle

\maketitle

\vskip-1cm

\begin{quote}

{\bf Abstract.}  We obtain a complete description of anisotropic scaling limits and the existence of scaling transition for a class
of negatively dependent linear random fields $\mathfrak{X}$
on ${\mathbb{Z}}^2$ with moving-average coefficients $a(t,s)$
decaying as $|t|^{-q_1}$ and $|s|^{-q_2}$ in the horizontal and vertical directions,  $q_1^{-1} + q_2^{-1} < 1 $ and satisfying
$\sum_{(t,s) \in \Z^2} a(t,s) = 0 $.
The scaling limits are taken over rectangles  whose sides
increase as $\lambda $ and $\lambda^\gamma $ when $\lambda \to \infty$, for any $\gamma >0$. The scaling  transition 
occurs
at $\gamma^\mathfrak{X}_0 >0$ if the scaling limits of $\mathfrak{X}$ are different and do not depend on $\gamma $ for $\gamma > \gamma^\mathfrak{X}_0 $ and $\gamma < \gamma^\mathfrak{X}_0$.
We prove that the scaling transition
in this model
is closely related to the presence
or absence of the edge effects. The paper extends the results in Pilipauskait\.e  and Surgailis (2017)
on the scaling transition for a related class of
random fields with long-range dependence.

\end{quote}

\section{Introduction}

A stationary random field (RF) $\mathfrak{X} = \{\mathfrak{X}({\mbt}); {\mbt} \in \Z^\nu\}$ on $\nu$-dimensional lattice $\Z^\nu, \nu \ge 1 $
with finite variance
is said to be (covariance) {\it long-range dependent} (LRD) if
$\sum_{\boldsymbol{t} \in \Z^\nu}  |{\rm Cov} (\mathfrak{X}({\mbt}), \mathfrak{X}({\mbf 0}))| = \infty  $, \
{\it short-range dependent} (SRD) if \ $\sum_{{\mbt} \in \Z^\nu}  |{\rm Cov} (\mathfrak{X}({\mbt}), \mathfrak{X}({\mbf 0}))| < \infty, \
\sum_{{\mbt} \in \Z^\nu}  {\rm Cov} (\mathfrak{X}({\mbt}), \mathfrak{X}({\mbf 0})) \ne 0$, and
{\it negatively dependent} (ND) if \ $\sum_{{\mbt} \in \Z^\nu}  |{\rm Cov} (\mathfrak{X}({\mbt}), \mathfrak{X}({\mbf 0}))| < \infty, \
\sum_{{\mbt} \in \Z^\nu}  {\rm Cov} (\mathfrak{X}({\mbt}), \mathfrak{X}({\mbf 0})) = 0$. The above definitions apply {\it per se} to RFs
with finite 2nd moment; related albeit not equivalent definitions of LRD, SRD, and ND properties
are discussed in \cite{book2012, lah2016,
ps2016} and other works.
For linear (moving-average) RFs,  Lahiri and Robinson \cite{lah2016} define similar concepts through summability properties
of the moving-average coefficients. The last paper also discusses the importance
of spatial LRD in applied sciences, including  the relevant literature.

The above classification plays an important role in limit theorems. Consider
the  sum $S^\mathfrak{X}_{K_\lambda}:= \sum_{{\mbt} \in K_\lambda} \mathfrak{X}({\mbt})$ of the values of RF $\mathfrak{X}$
over large `sampling region'  $K_\lambda \subset \Z^\nu $ with
$|K_\lambda| = \sum_{{\mbt} \in K_\lambda } 1 \to \infty \ (\lambda \to \infty) $. Under additional conditions,
the variance ${\rm Var}(S^\mathfrak{X}_{K_\lambda})$ grows faster than
$|K_\lambda|$ under LRD, as $ O(|K_\lambda|) $ under SRD, and slower than $|K_\lambda|$ under ND. In the latter case,
${\rm Var}(S^\mathfrak{X}_{K_\lambda})$
may grow as slow as the `volume' $|\partial K_\lambda| $ of the boundary
$\partial K_\lambda $,
or even slower than  $|\partial K_\lambda| $, giving rise to  `edge effects' which may
affect or dominate
the limit distribution of $S^\mathfrak{X}_{K_\lambda}$; see \cite{lah2016, rip1988}.

Probably, the most studied case of limit theorems for RFs deal with rectangular summation regions, which allows
for partial sums and limit RFs, similarly as in the case $\nu=1$.
Let $\mathfrak{X} = \{\mathfrak{X}(\mbt); \mbt \in \Z^\nu \}$ be a stationary random  field (RF) on $\Z^\nu, \, \nu \ge 1, \,
\mbgamma = (\gamma_1, \cdots, \gamma_\nu) \in \R^\nu_+ $ be a collection of positive numbers,
and
\begin{eqnarray} \label{Kp}
K_{\la, \smbgamma} (\mbx)&:=&[1, \lfloor \la^{\gamma_1} x_1\rfloor ] \times \cdots \times [1, \lfloor \la^{\gamma_\nu} x_\nu\rfloor ], \qquad  \mbx = (x_1, \cdots, x_\nu) \in \R^\nu_+,
\end{eqnarray}
be a family of $\nu$-dimensional  `rectangles' indexed by $\la >0$,
whose sides grow at generally different rate $O(\la^{\gamma_i}), i=1,\cdots, \nu$ as $\la \to \infty$, and
\begin{eqnarray}\label{Sn}
S^\mathfrak{X}_{\la, \mbgamma}(\mbx)
&:=&\sum_{{\mbt} \in K_{\la, \mbgamma} ({\mbx})  }  \mathfrak{X}(\mbt), \qquad \mbx \in \R^\nu_+
\end{eqnarray}
be the corresponding partial  sums RF.
\cite{ps2015, ps2016, pils2017, sur2019}  discussed the anisotropic  scaling limits  for any
$\mbgamma \in \R^\nu_+  $
of some classes of
LRD RFs $\cY$ in dimension $\nu = 2,3$, viz.,
\begin{equation}\label{partS}
A^{-1}_{\la, \mbgamma} S^\mathfrak{X}_{\la,\mbgamma}(\mbx)  \ \limfdd \ V^\mathfrak{X}_{\mbgamma} (\mbx), \quad \mbx \in \R^\nu_+
\end{equation}
as $\la  \to \infty$,  where $A_{\la, \mbgamma} \to  \infty $ is a normalization.
Following \cite{pils2016, sur2019}  the family
$\{ V^\mathfrak{X}_{\mbgamma};  \mbgamma \in \R^\nu_+\} $ of all scaling limits in \eqref{partS} will
be called the {\it scaling diagram of RF $\mathfrak{X}$}. \cite{pils2016} noted that the scaling diagram provides a more complete
`large-scale summary' of RF $\mathfrak{X}$ compared to (usual)  isotropic or anisotropic scaling at fixed $\mbgamma \in \R^\nu_+$
as discussed in \cite{AnhLRM2012, dobmaj1979,
lah2016, lav2007, leo1999, sur1982} and other works.

\cite{ps2015, ps2016, pils2017} observed that for a large class of LRD RFs  $\mathfrak{X}$ in dimension $\nu =2$, the scaling
diagram essentially consists of three elements:
$V^\mathfrak{X} = \{V^\mathfrak{X}_0, V^\mathfrak{X}_+, V^\mathfrak{X}_- \}$,  $V^\mathfrak{X}_0$ termed the {\it well-balanced} and
$V^\mathfrak{X}_\pm $ the {\it unbalanced} scaling limits of $\mathfrak{X}$. For $\nu=2$ without loss of generality (w.l.g.)
we can assume $\mbgamma = (1, \gamma) \in \R^2_+ $ or $\gamma_1 = 1, \gamma_2 = \gamma $   in \eqref{Sn}
and denote
\begin{eqnarray} \label{SXdef}
&S^\mathfrak{X}_{\la,\gamma}(x,y) := \sum_{(t,s) \in K_{[\la x, \la^\gamma y]}} \mathfrak{X}(t,s), \qquad (x,y) \in \R^2_+
\end{eqnarray}
and $V^\mathfrak{X}_\gamma, A_{\la, \gamma} $  the corresponding sum in \eqref{Sn} and  
quantities in \eqref{partS} defined
for $\gamma >0$, where
$K_{[\la x, \la^\gamma y]} = \{ (t,s) \in \Z^2: 1 \le t \le \lfloor \la x\rfloor,
1 \le s \le \lfloor \la^\gamma y \rfloor \} $.
\cite{ps2015, ps2016, pils2017} proved
that there exists a (nonrandom) $\gamma^\mathfrak{X}_0 >0$ such that $V^\mathfrak{X}_\gamma$ do not depend on $\gamma $
for  $\gamma > \gamma^\mathfrak{X}_0$ and
$\gamma < \gamma^\mathfrak{X}_0$, viz.,
\begin{equation}\label{VX0}
V^\mathfrak{X}_\gamma = \begin{cases} V^\mathfrak{X}_+,  &\gamma > \gamma^\mathfrak{X}_0, \\
V^\mathfrak{X}_-, &\gamma < \gamma^\mathfrak{X}_0, \\
V^\mathfrak{X}_0, &\gamma = \gamma^\mathfrak{X}_0
\end{cases}
\end{equation}
and $V^\mathfrak{X}_+ \neqfdd a V^\mathfrak{X}_- \ (\forall \, a >0)$.
The above fact was termed the {\it scaling transition} \cite{ps2015, ps2016}.
It was noted in the above-mentioned works that the scaling transition constitutes a new and general feature of spatial
dependence which occurs in many spatio-temporal models including telecommunications and economics
\cite{gaig2003, kajt2008, miko2002, pils2014, pils2015,  pils2016, lei2018, pilss2020}.
However, as noted in \cite{pils2017, sur2019},
these studies were limited to LRD models and
the existence of the scaling transition under ND remained open.

The present paper discusses the scaling transition for linear ND RFs on $\Z^2 $ having a moving-average
representation
\begin{equation}\label{Xlin}
\mathfrak{X}(t,s) \ = \ \sum_{(u,v) \in \Z^2} a(t-u, s-v) \vep(u,v), \qquad (t,s) \in \Z^2,
\end{equation}
in  standardized i.i.d.\ sequence $\{ \vep(u,v); (u,v) \in \Z^2\}, \E \vep (u,v) =0, \E \vep^2(u,v) =1$ with deterministic
moving-average coefficients
\begin{equation}\label{acoefL}
a(t,s)\  = \  \frac{1}{ (|t|^2 +  |s|^{2q_2/q_1})^{q_1/2}}\Big(L_0 \big(\frac{t}{(|t|^{2} +|s|^{2q_2/q_1})^{1/2}}\big) + o(1)\Big),
\qquad |t|+|s| \to  \infty,
\end{equation}
$(t,s) \neq (0,0), $ where $q_i>0, i=1,2 $ satisfy
\begin{equation}\label{qLRD}
0 \ < \ Q := \frac{1}{q_1} + \frac{1}{ q_2} \ < \ 2
\end{equation}
which imply $q_i > 1/2, i=1,2$. In  \eqref{acoefL},
$L_0(u), u \in [-1,1]$ is a bounded piece-wise continuous function on $[-1,1]$ termed the  {\it angular function} in \cite{pils2017}.
(We note that the boundedness and continuity assumptions on the angular function
do not seem necessary for our results and possibly can be relaxed.)
The form of moving-average coefficients in \eqref{acoefL} is the same as in \cite{pils2017} and
can be generalized to some extent but we prefer to  use  \eqref{acoefL}
for better comparison with the results of \cite{pils2017}. Condition $Q< 1 $ guarantees that
$\sum_{(t,s)\in \Z^2} |a(t,s)| < \infty $ and the ND property of $\mathfrak{X}$ in \eqref{Xlin} is a consequence of the zero-sum condition:
\begin{eqnarray}\label{zero}
\sum_{(t,s) \in \Z^2} a(t,s)&=&0, \qquad \text{for} \quad Q < 1.
\end{eqnarray}
In contrast,
\cite{pils2017} assumes $1 < Q < 2 $ implying $\sum_{(t,s)\in \Z^2} |a(t,s)| = \infty $ and the LRD property of the corresponding
linear RF $\mathfrak{X}$ in \eqref{Xlin}. The linear model in \eqref{Xlin} and the results of this paper can be regarded as an extension of the classical
results of Davydov \cite{dav1970} in the time series setting, who identified all partial sums limits (fractional Brownian motions)
of moving-average processes with one-dimensional `time'.

The main results of this paper and \cite{pils2017} are illustrated in Fig. 1 showing 8 regions $R_{11}, \cdots, R_{33} $ of the parameter set
$\{(1/q_1, 1/q_2): 0 < Q < 2 \} $ of the linear RF $\mathfrak{X}$ in \eqref{Xlin}-\eqref{zero}
 with different unbalanced limits. We remark that when
$\sum_{(t,s) \in \Z^2} a(t,s)\ne 0 $ and $Q< 1 $ the linear RF  $\mathfrak{X}$ in \eqref{Xlin}-\eqref{acoefL} is SRD and all scaling limits $V^\mathfrak{X}_\gamma, \gamma >0$
agree with Brownian sheet $B_{1/2,1/2} $, see (\cite{pils2017}, Theorem~3.4),  meaning that in this case the  scaling transition does not occur.

\begin{center}
\begin{figure}[ht]
\begin{center}
\includegraphics[width=15 cm,height=17cm]{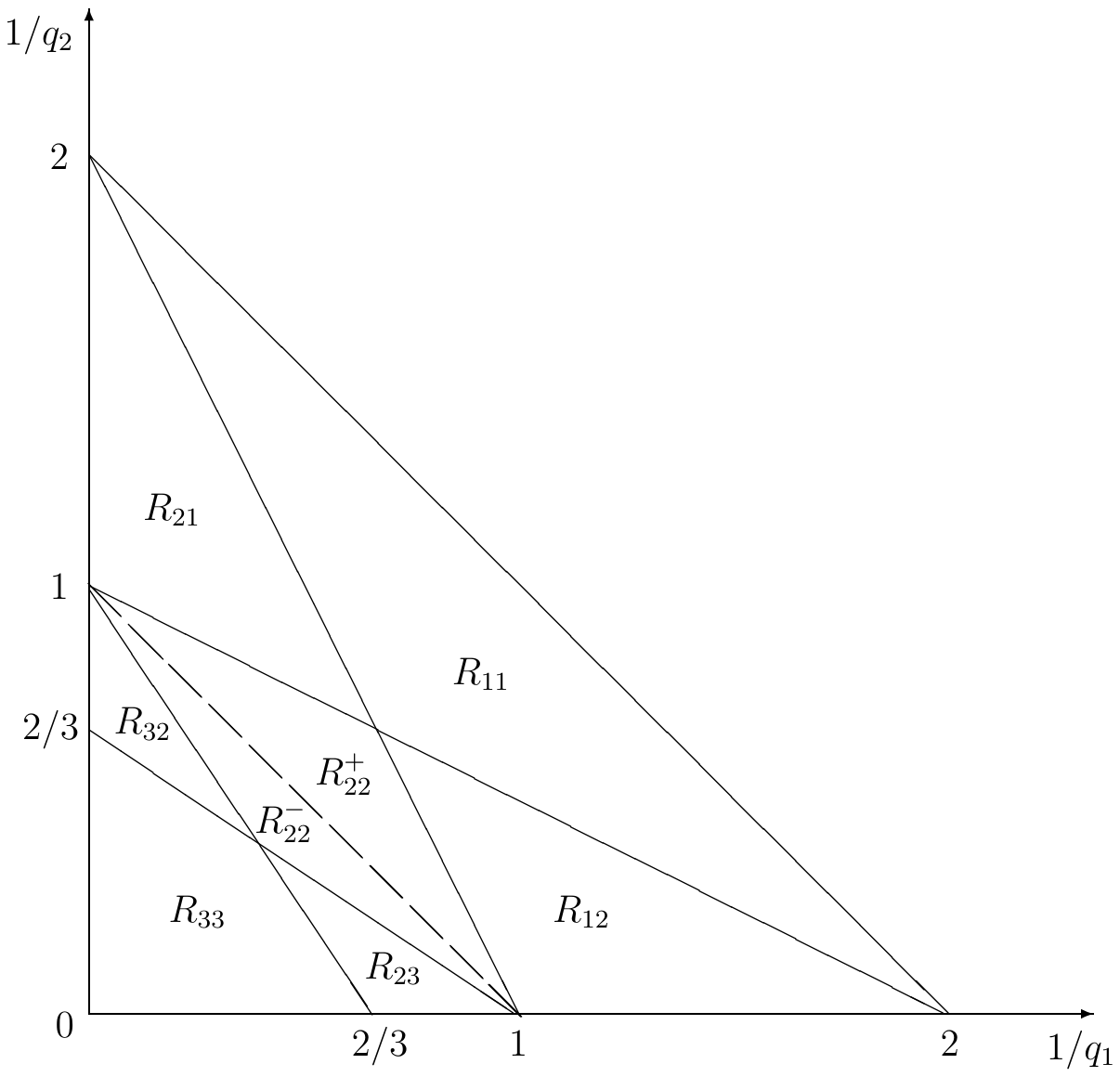}
\vspace{-9.5cm}
\end{center}
\end{figure}
\end{center}

\vskip-1cm

\hskip2cm Figure 1. \parbox[t]{12cm} {\small Regions in the parameter set
$0 < Q < 2 $ with different unbalanced limits of RF $\mathfrak{X}$ in \eqref{Xlin}-\eqref{zero}.
 The dashed segment separates the LRD and ND regions.}

\bigskip

The regions $R_{11}, \cdots, R_{33}$ in Figure~1 are described in Table~1. Recall the definition
of fractional Brownian sheet (FBS) $B_{H_1,H_2} = \{B_{H_1,H_2}(x,y); (x,y) \in \R^2_+ \} $ with Hurst parameters
$0< H_1, H_2 \le 1 $ as a Gaussian process with zero mean and covariance
\begin{equation}\label{FBScov}
\E B_{H_1,H_2}(x_1,y_1) B_{H_1,H_2}(x_2,y_2) = (1/4) (x_1^{2H_1} + x_2^{2H_1} - |x_1-x_2|^{2H_1})
 (y_1^{2H_2} + y_2^{2H_2} - |y_1-y_2|^{2H_2}),
\end{equation}
see \cite{aya2002}. FBS $B_{H_1,H_2}$ with one of the parameters $H_i$ equal to $1/2$ or $1$ have a very specific
dependence structure (either independent or completely dependent (invariant) increments in one direction, see
\cite{ps2015, ps2016})  and play a particular role in our work. We extend the above class of Gaussian RFs to parameter values $(H_1,H_1) \in [0,1]^2 $
by setting
\begin{eqnarray} \label{BH0}
\E B_{H_1,0}(x_1,y_1) B_{H_1,0}(x_2,y_2)
&:=&\lim_{H_2 \to +0}  \E B_{H_1,H_2}(x_1,y_1) B_{H_1,H_2}(x_2,y_2),   \quad 0< H_1 \le 1, \\
\E B_{0,H_2}(x_1,y_1) B_{0,H_2}(x_2,y_2)
&:=&\lim_{H_1 \to +0}  \E B_{H_1,H_2}(x_1,y_1) B_{H_1,H_2}(x_2,y_2),   \quad 0< H_2 \le 1, \nn \\
\E B_{0,0}(x_1,y_1) B_{0,0}(x_2,y_2)
&:=&\lim_{H_1, H_2 \to +0}  \E B_{H_1,H_2}(x_1,y_1) B_{H_1,H_2}(x_2,y_2),  \nn
\end{eqnarray}
for any $(x_i,y_i) \in \R^2_+, i=1,2$. The limit covariances on the r.h.s. of \eqref{BH0} are explicitly written in Remark \ref{rem0}.
FBS $B_{H_1,H_2}$ where one or both parameters are equal to $0$ are rather unusual (have nonseparable paths).
The appearance of such RFs in limit theorems is surprising but seems to be rather common under edge effects,
see Remark \ref{remB0}.


\begin{table}[htbp]

\begin{center}
\begin{tabular}{| c | c | c | c | c |}
\hline
{Parameter region} & {Critical $\gamma^\mathfrak{X}_0$} & {$V^\mathfrak{X}_+$} & {$ V^\mathfrak{X}_-$}  & {Range of Hurst parameters}   \\
\hline
{$R_{11}$} & $q_1/q_2$ & $B_{1,\tilde H_2}$  & $B_{\tilde H_1,1}$ & $1/2 < \tilde H_1, \tilde H_2 < 1$ \\
\hline
$R_{12}$ & $q_1/q_2$ & $B_{H_1, 1/2}$   &  $B_{\tilde H_1,1}$ &   $1/2 < \tilde H_1, H_1 < 1$   \\
\hline
$R_{21}$ & $q_1/q_2$ & $B_{1, \tilde H_2}$   &  $B_{1/2, H_2}$ &   $1/2 < \tilde H_2, H_2 < 1$   \\
\hline
$R^+_{22}$ & $q_1/q_2$ & $B_{H_1, 1/2}$   &  $B_{1/2, H_2}$ &   $1/2 < H_1, H_2 < 1$   \\
\hline
$R^-_{22}$ & $q_1/q_2$ & $B_{H_1, 1/2}$   &  $B_{1/2, H_2}$ &   $0 < H_1, H_2 < 1/2$   \\
\hline
$R_{23}$ & $2q_1(1-Q)$ & $B_{H_1, 1/2}$   &  $B_{1/2,0}$ &   $0 < H_1 < 1/2$   \\
\hline
$R_{32}$ & $1/2q_2(1-Q)$ & $B_{0,1/2}$   &  $B_{1/2, H_2}$ &   $0 < H_2 < 1/2$   \\
\hline
$R_{33}$ & $1$ & $B_{0, 1/2}$   &  $B_{1/2,0}$ &     \\
\hline
\end{tabular} \\
\caption{Unbalanced scaling limits $V^\mathfrak{X}_\pm $ in regions $R_{11}, \cdots, R_{33}$ in Fig. 1
 }\label{tab:1}
\end{center}

\end{table}

Parameters $\tilde H_i, H_i, i=1,2 $ in Table~1 (expressed in terms of $q_1, q_2$)
are specified in the beginning of Sec.~2.
The description in Table~1 is not very precise since it omits various asymptotic constants which may vanish in some cases,
meaning that some additional conditions on $a(t,s)$ are needed for the validity of the results in this table.
The rigorous formulations including the normalizing constants $A_{\lambda, \gamma} $ are presented in Sec.~2.
The limit distributions in regions $R_{11}, R_{12}, R_{21}, R^+_{22}$ refer to LRD set-up
and are part of \cite{pils2017}. The new results under ND  refer to regions $R_{22}^-, R_{23}, R_{32}, R_{33}$
of Table~1. Note $2q_1 (1-Q) > q_1/q_2 $ in $R_{23} $ and  $1/(2q_2 (1-Q)) < q_1/q_2 $ in  $R_{32} $. The limit $B_{1/2,0}$ in
region $R_{23} $ can be related to the
`horizontal edge effect' which dominates the limit distribution of $S^\mathfrak{X}_{\lambda, \gamma} $ of \eqref{SXdef}
unless the vertical length $O(\lambda^\gamma)$ of
$K_{[\lambda x, \lambda^\gamma y]} $
grows fast enough vs. its horizontal length $O(\lambda)$, or
$\gamma  > 2q_1 (1-Q)$ holds, in which case  FBS $B_{H_1, 1/2}$ dominates.  
Similarly, $B_{0, 1/2}$ in
region  $R_{32} $ can be related to the
`vertical edge effect' appearing in the limit of $S^\mathfrak{X}_{\lambda, \gamma} $
unless the vertical length $O(\lambda^\gamma)$ increases sufficiently slow w.r.t the horizontal length $O(\lambda)$, or
$\gamma  < 1/2q_2 (1-Q)$ holds, in which case FBS $B_{1/2, H_2}$ dominates.
Finally, $R_{33}$ can be characterized as the parameter region
where the edge effects (either horizontal, or vertical) completely dominate the limit behavior of $S^\mathfrak{X}_{\lambda, \gamma} $.
The above interpretation of $R_{32}, R_{23} $ and $R_{33}$  is based on the approximations of  $S^\mathfrak{X}_{\lambda, \gamma} $
by suitable `edge terms' which are discussed in Sec.~3.

The results of the present work are related to the work Lahiri and Robinson \cite{lah2016} which
discussed the limit distribution of sums of linear LRD, SRD and ND RFs over homothetically inflated or isotropically  rescaled (i.e., $\gamma_1 = \cdots = \gamma_\nu = 1$) star-like  regions $K_\lambda $ of very general form. This generality of  $K_\lambda $ does not seem to allow for a natural
introducing of partial sums, restricting the problem to the convergence of one-dimensional distributions in contrast to finite-dimensional
distributions in the present paper.
While  \cite{lah2016} consider several forms of moving-average coefficients,
the only case when $a(t,s)$ in \eqref{acoefL} satisfy the assumptions in \cite{lah2016} seems to be the `isotropic' case $q_1 = q_2$. As explained
in Remark \ref{LRob},  in the case $q_1 = q_2 $ and $\mbgamma = (1,1) $
our limit results agree with  \cite{lah2016}, including the `edge effect'.
We also note Damarackas and Paulauskas \cite{dam2017} who discussed partial sums limits of linear LRD, SRD and ND RFs, possibly
with infinite variance and  moving-average coefficients which factorize along coordinate axes (i.e., different from
\eqref{acoefL}) in which case the scaling limits in \eqref{partS} do not depend on $\mbgamma $ and the scaling transition
does not exist. See also (\cite{sur2019}, Remark~4.1).

The  rest of the paper is organized as follows. Sec.~2 contains the main results (Theorems 2.1 - 2.3). The proofs
of these results are given in Sec.~3. Sec.~4 presents two examples of fractionally integrated ND
RFs, extending the examples of fractionally integrated LRD RFs in \cite{pils2017}.
Sec.~5 (Concluding remarks) summarizes the  main discoveries and mentions several open problems.

\smallskip

 {\it Notation.} In what follows, $C$ denote  generic positive constants
which may be different at different locations. We write $\limfdd,   \eqfdd, $ and $ \neqfdd $
for the weak convergence, equality and inequality
of finite-dimensional
distributions, respectively.
$\R^\nu_+ := \{  \mbx =  (x_1, \cdots, x_\nu) \in \R^\nu: x_i > 0,
i=1, \cdots, \nu \},
\R_+ := \R^1_+$,
$\lfloor x \rfloor := \max \{ k \in \Z: k \le x \},  \lceil x \rceil := \min \{ k \in \Z: k \ge x \}, x \in \R$.
$\1(A)$ stands for the  indicator function of a set $A$.

\section{Main results}

Throughout the paper we use the following notation:
\begin{eqnarray}\label{Hnota}
&&H_1:=  \frac{3}{2} + \frac{q_1}{q_2}  - q_1 \  = \  \frac{1}{2} + q_1 (Q-1), \qquad
H_2:= \frac{3}{2} + \frac{q_2}{q_1}  - q_2 \  = \  \frac{1}{2} + q_2 (Q-1),  \\
&&\tilde H_1 := \frac{3}{2} - q_1 + \frac{q_1}{2q_2}  \ = \  1 - \frac{q_1}{2} (2-Q), \qquad
\tilde H_2 := \frac{3}{2} - q_2 + \frac{q_2}{2q_1}  \ = \  1 - \frac{q_2}{2} (2-Q), \nn \\
&&\gamma^0:=\frac{q_1}{q_2}, \quad
\gamma^0_{\rm edge,1}:= \frac{1}{2q_2(1-Q)}, \quad \gamma^0_{\rm edge,2} := 2q_1(1-Q), \nn \\
&&Q_{\rm edge,1}:= \frac{3}{2q_1} + \frac{1}{q_2},  \qquad  Q_{\rm edge,2}:= \frac{1}{q_1} + \frac{3}{2q_2},
\qquad \widetilde Q_1 := \frac{1}{2q_1} + \frac{1}{q_2}, \qquad
\widetilde Q_2 := \frac{1}{q_1} + \frac{1}{2q_2}.
 \nn
\end{eqnarray}
We recall that $q_i, i=1,2 $ are the asymptotic parameters of the MA coefficients in \eqref{acoefL} satisfying \eqref{zero}
and $Q $ is defined in \eqref{qLRD}.
$H_i, \tilde H_i, i=1,2 $ in \eqref{Hnota} are the Hurst parameters of  the limit FBS in Table~1 and in the subsequent Theorems
\ref{main0}--\ref{main2}. The angular function $L_0$ in \eqref{acoefL} enters the expression of the kernel
$a_\infty $ in \eqref{a0} through which these limit FBS are defined.
 The positive quantities  $\gamma^0, \gamma^0_{{\rm edge},i}, i=1,2 $ correspond to the critical points
$\gamma^\mathfrak{X}_0$ in \eqref{VX0}. Also
note that the partition in Figure~1 is formed by segments belonging to the lines $Q_{{\rm edge},i} =1, \widetilde Q_i = 1, i=1,2$ in \eqref{Hnota} and
$Q = 1, Q=2$.  Let
\begin{eqnarray} \label{a0}
a_\infty (t,s) &:=&  (|t|^2 +  |s|^{2q_2/q_1})^{-q_1/2} L_0 \big(t/(|t|^{2} +|s|^{2q_2/q_1})^{1/2}\big), \qquad (t,s) \in \R^2.
\end{eqnarray}
Note $a_\infty $ is the scaling limit of $a(t,s)$ in \eqref{acoefL}: 
$\lim_{\lambda \to \infty} \lambda^{q_1} a(\lceil \lambda t\rceil, \lceil \lambda^{q_1/q_2} s\rceil ) $  $= a_\infty (t,s)$
for any $(t,s) \in \R^2, (t,s) \ne (0,0).$
We also use the following notation for kernels of the limit RFs expressed as stochastic integrals w.r.t. Gaussian white noise
$W(\d u, \d v), (u,v) \in \R^2,  \E W(\d u, \d v) = 0,  \E (W(\d u, \d v))^2 = \d u \d v $
on $\R^2$. Namely,
for $(u,v) \in \R^2, (x,y) \in \R^2_+$,  let
\begin{eqnarray} \label{h00}
h_{0}(x,y; u,v)&:=&\begin{cases}
\int_{(0,x]\times (0,y]} a_\infty (t-u,s-v) \d t \d s,
&(u,v) \not\in (0,x] \times (0,y], \\
-\int_{\R^2 \setminus (0,x]\times (0,y]} a_\infty (t-u,s-v) \d t \d s,
&(u,v) \in (0,x] \times (0,y],
\end{cases} \\
h_{1}(x,y; u,v)&:=&\1(0< v \le y)
\begin{cases}
\int_{(0,x]\times \R} a_\infty (t-u,s) \d t \d s,
&u \not\in (0,x], \\
-\int_{(\R \setminus (0,x])\times \R} a_\infty (t-u,s) \d t \d s,
&u \in (0,x],
\end{cases} \nn \\
h_{2}(x,y; u,v)&:=&\1(0< u \le x)
\begin{cases}
\int_{\R \times (0,y]} a_\infty (t,s-v) \d t \d s,
&v \not\in (0,y], \\
-\int_{\R \times (\R \setminus (0,y])} a_\infty (t,s-v) \d t \d s,
&v \in (0,y].
\end{cases} \nn
\end{eqnarray}
We point out that definitions {\it \eqref{h00}  apply under ND condition $Q < 1 $} only.
Under the LRD condition $1 < Q < 2 $ the corresponding kernels take a somewhat simpler form:
\begin{eqnarray} \label{h00LRD}
h_{0}(x,y; u,v)&:=&
\int_{(0,x]\times (0,y]} a_\infty (t-u,s-v) \d t \d s,  \\
h_{1}(x,y; u,v)&:=&\1(0< v \le y)
\int_{(0,x]\times \R} a_\infty (t-u,s) \d t \d s, \nn \\
h_{2}(x,y; u,v)&:=&\1(0< u \le x)
\int_{\R \times (0,y]} a_\infty (t,s-v) \d t \d s, \nn
\end{eqnarray}
\eqref{h00LRD}  hold for any $(u,v) \in \R^2 $. In \eqref{hFBM} the kernels $h_i, i=1,2 $ of \eqref{h00} 
are explicitly written. For the sake of completeness, we also introduce the kernels
\begin{eqnarray} \label{hLRD}
&&\tilde h_{1}(x,y; u,v):= x
\int_{(0,y]} a_\infty (u,s-v) \d s,  \qquad
\tilde h_{2}(x,y; u,v):=  y \int_{(0,x]} a_\infty (t-u,v) \d t
\end{eqnarray}
used in the definition of the unbalanced limits in regions $R_{11}, R_{12}, R_{21}$ of Table~1. Using the above notation
we define the well-balanced limit
\begin{eqnarray}\label{V0}
V^\mathfrak{X}_0 (x,y)&:=&\int_{\R^2} h_0(x,y; u,v) W(\d u, \d v), \quad \text{for}  \quad   Q_{\rm edge,1} \wedge Q_{\rm edge, 2} > 1,  0 < Q < 2, Q\ne 1
\end{eqnarray}
in both cases $0 < Q < 1$ (ND) and $1< Q < 2 $ (LRD), with respective $h_0$ \eqref{h00} and \eqref{h00LRD}.
Similarly, we define RFs
\begin{eqnarray} \label{V+}
V^\mathfrak{X}_+ (x,y)&:=&\int_{\R^2} h_1(x,y; u,v) W(\d u, \d v), \quad \text{for}  \
\left\{\begin{array}{ll}
Q < 1, &Q_{\rm edge,1} > 1, \\
1 < Q < 2, &\widetilde Q_{1} < 1
\end{array}\right.
\end{eqnarray}
and
\begin{eqnarray} \label{V-}
V^\mathfrak{X}_- (x,y)&:=&\int_{\R^2} h_2(x,y; u,v) W(\d u, \d v), \quad \text{for}  \
\left\{\begin{array}{ll}
Q < 1, &Q_{\rm edge,2} > 1, \\
1 < Q < 2, &\widetilde Q_{2} < 1.
\end{array}\right.
\end{eqnarray}
Finally, we define RFs
\begin{eqnarray} \label{tiV}
V^\mathfrak{X}_+ (x,y)&:=&\int_{\R^2} \tilde h_1(x,y; u,v) W(\d u, \d v), \quad \text{for}  \
1 < Q < 2, \ \widetilde Q_{1} > 1, \\
V^\mathfrak{X}_- (x,y)&:=&\int_{\R^2} \tilde h_2(x,y; u,v) W(\d u, \d v), \quad \text{for}  \
1 < Q < 2, \  \widetilde Q_{2} > 1. \nn
\end{eqnarray}
Recall that the stochastic integral $I(h) := \int_{\R^2} h(u,v) W(\d u, \d v) $ w.r.t. the white noise $W$  is well-defined
for any $h \in L^2(\R^2)$ and has a Gaussian distribution with zero mean and variance
$\E I(h)^2 = \|h\|^2 := \int_{\R^2} h(u,v)^2 \d u \d v $.

\begin{prop} \label{Vex} The Gaussian RFs  $V^\mathfrak{X}_0  = \{ V^\mathfrak{X}_0 (x,y); (x,y) \in \R^2_+\},
V^\mathfrak{X}_\pm  = \{ V^\mathfrak{X}_\pm (x,y); (x,y) \in \R^2_+\} $ are well-defined in respective regions of parameters $q_i >0, i=1,2$ as indicated
in \eqref{V0},  \eqref{V+},  \eqref{V-} and \eqref{tiV}. Moreover,
\begin{eqnarray}\label{VB1}
V^\mathfrak{X}_+ &\eqfdd&\begin{cases}
\sigma_1 B_{H_1,1/2}, &Q_{\rm edge, 1}>1, \widetilde Q_1 < 1,  Q \ne 1,\\
\tilde \sigma_1 B_{1, \tilde H_2}, &\widetilde Q_1 > 1
\end{cases}
\end{eqnarray}
and
\begin{eqnarray} \label{VB2}
V^\mathfrak{X}_- &\eqfdd&\begin{cases}
\sigma_2 B_{1/2, H_2}, &Q_{\rm edge, 2}>1,  \widetilde Q_2 < 1,  Q \ne 1,\\
\tilde \sigma_2 B_{\tilde H_1,1}, &\widetilde Q_2 > 1
\end{cases}
\end{eqnarray}
where $\sigma_i := \|h_i(1,1;\cdot,\cdot)\| < \infty, \tilde \sigma_i := \|\tilde h_i(1,1;\cdot,\cdot)\| < \infty, i=1,2. $

\end{prop}

We note that for $1< Q < 2 $ the statement of Proposition \ref{Vex} is part of (\cite{pils2017}, Theorems~3.1-3.3).
For $Q < 1 $ and $h_0 $ in \eqref{h00} the fact that  $h_0 \in L^2 (\R^2) $ (i.e., that $V^\mathfrak{X}_0$ is well-defined)
follows from Lemma \ref{Lem}, \eqref{aint8}-\eqref{aint9}. The statements \eqref{VB1} and \eqref{VB2} for $Q < 1 $ follow by
rewriting the kernels $h_i, i=1,2$ in \eqref{h00} as
\begin{eqnarray*}
h_1(x,y;u,v)&=&\1(0< v \le y)
\begin{cases}
\int_{(0,x]} |t-u|^{-q_1(1- \frac{1}{q_2})} L_{1,{\rm sign}(t-u)} \d t, &u \not\in (0,x], \\
-\int_{\R\setminus (0,x]} |t-u|^{-q_1(1- \frac{1}{q_2})} L_{1,{\rm sign}(t-u)} \d t, &u \in (0,x],
\end{cases}, \\  
h_2(x,y;u,v)&=&\1(0< u  \le x)
\begin{cases}
\int_{(0,y]} |s-v|^{-q_2(1- \frac{1}{q_1})} L_{2, {\rm sign}(s-v)}
\d s, &v \not\in (0,y], \\
- \int_{\R\setminus (0,y]} |s-v|^{-q_2(1- \frac{1}{q_1})} L_{2, {\rm sign}(s-v)}
\d s, &v \in (0,y],
\end{cases}  
\end{eqnarray*}
where
\begin{equation*}
L_{1, \pm} := \int_{\R} a_\infty (\pm 1, s) \d s, \qquad  L_{2,\pm} := \int_{\R} a_\infty (t,\pm 1) \d t
\end{equation*}
(note $L_{2,+} = L_{2,-}$ since  $a_\infty (t,s) = a_\infty (t,-s)$  is symmetric in $s$).
Furthermore, the above integrals can be rewritten as
\begin{eqnarray} \label{hFBM}
h_1(x,y;u,v)&=&\frac{\1(0< v \le y)}{(1/2) - H_1}
\Big\{ L_{1,+} \big((x-u)_+^{H_1-(1/2)} - (-u)_+^{H_1- (1/2)}\big) \\
&&\hskip2cm + \ L_{1,-} \big((x-u)_-^{H_1- (1/2)} - (-u)_-^{H_1-(1/2)}\big)\Big\},  \nn \\
h_2(x,y;u,v)&=&\frac{\1(0< u \le x)}{(1/2) - H_2}
\Big\{ L_{2,+} \big((y-v)_+^{H_2-(1/2)} - (-v)_+^{H_2-(1/2)}\big) \nn \\
&&\hskip2cm + \ L_{2,-} \big((y-v)_-^{H_2-(1/2)} - (-v)_-^{H_2-(1/2)}\big)\Big\}, \nn
\end{eqnarray}
with $H_i, i=1,2 $ defined in \eqref{Hnota}. Whence we immediately see that $h_i, i=1,2$ factorize into a product
of kernels of (fractional) Brownian motions with one-dimensional time, see \cite{sam1994},  implying facts
\eqref{VB1} and \eqref{VB2} for $Q < 1 $.

The main object of this paper is  establishing  the scaling transition for ND RF $\mathfrak{X}$, viz.,
\begin{equation}\label{partS1}
A^{-1}_{\la, \gamma} S^\mathfrak{X}_{\la,\gamma}(x,y)  \ \limfdd \
\begin{cases} V^\mathfrak{X}_+(x,y),  &\gamma > \gamma^\mathfrak{X}_0, \\
V^\mathfrak{X}_-(x,y), &\gamma < \gamma^\mathfrak{X}_0, \\
V^\mathfrak{X}_0(x,y), &\gamma = \gamma^\mathfrak{X}_0,
\end{cases}
\end{equation}
where $A_{\lambda,\gamma} \to \infty \, (\lambda \to \infty)$ is normalization and  $S^\mathfrak{X}_{\la,\gamma}(x,y) $ is defined in
\eqref{SXdef}.

This question is treated in Theorems \ref{main0}-\ref{main2}. In these theorems, $\mathfrak{X}$ is a linear RF as in \eqref{Xlin} with
moving-average coefficients $a(t,s)$ satisfying \eqref{acoefL}, \eqref{zero}, where the parameters  $q_1, q_2$ satisfy
$0 < Q < 1 $ and
belong to different regions  of Figure~1.

\begin{thm} \label{main0} {\rm [Region $R^-_{22}$:
`no edge effects'.]}
Let $Q_{\rm edge,1} \wedge Q_{\rm edge,2} >1. $
Then the convergence in \eqref{partS1} holds with $\gamma^\mathfrak{X}_0 = \gamma^0 = q_1/q_2,
A_{\la, \gamma} = \lambda^{H(\gamma)} $ and the limit RFs specified in \eqref{V+}, \eqref{V-},
viz.,
\begin{eqnarray*}
\begin{cases} V^\mathfrak{X}_+ = \sigma_1 B_{H_1,1/2}, \\
V^\mathfrak{X}_- = \sigma_2 B_{1/2, H_2},  \\
V^\mathfrak{X}_0 = V_0,
\end{cases} \qquad
H(\gamma) = \begin{cases}
H_1 + (\gamma/2), &\gamma > \gamma^0, \\
\gamma H_2 + (1/2), &\gamma < \gamma^0, \\
H_1 + (\gamma^0/2) = \gamma^0 H_2 + (1/2), &\gamma = \gamma^0
\end{cases}
\end{eqnarray*}
where $H_i, i=1,2 $ and $\sigma_i, i=1,2 $ defined in \eqref{Hnota} and in Proposition \ref{Vex},
respectively.

\end{thm}

Define
\begin{eqnarray}\label{sedge}
\sigma^2_{\rm edge, 1}&:=&2\sum_{v \ge 0} \big(\sum_{t\in \Z, s \ge 1} a(t,s+v) \big)^2 +
2\sum_{v \le -1} \big(\sum_{t\in \Z, s \le 0} a(t,s+v) \big)^2,  \\
\sigma^2_{\rm edge,2}&:=&2\sum_{u \ge 0} \big(\sum_{t\ge 1, s \in \Z} a(t+u,s) \big)^2 + 2\sum_{u \le -1} \big(\sum_{t\le 0, s \in \Z} a(t+u,s) \big)^2. \nn
\end{eqnarray}
The convergence of the above series in the corresponding regions of parameters $q_1, q_2$  will be established later.

\begin{thm} \label{main1} {\rm [Regions $R_{23}$ and $R_{32}$:
`one-sided edge effects'.]}

\noi (i)
Let $Q_{\rm edge,2} < 1  < Q_{\rm edge,1} $.
Then the convergence in \eqref{partS1} for $\gamma \ne \gamma^\mathfrak{X}_0$ holds
with $\gamma^\mathfrak{X}_0 = \gamma^0_{\rm edge,2},
A_{\la, \gamma} = \lambda^{H(\gamma)} $ and the limit RFs specified in \eqref{V+}, \eqref{V-},
viz.,
\begin{eqnarray*}
\begin{cases} V^\mathfrak{X}_+ = \sigma_1 B_{H_1,1/2}, \\
V^\mathfrak{X}_- = \sigma_{\rm edge,1} B_{1/2,0},
\end{cases} \qquad
H(\gamma) = \begin{cases}
H_1 + (\gamma/2), &\gamma > \gamma^0_{\rm edge,2}, \\
1/2, &\gamma < \gamma^0_{\rm edge,2},
\end{cases}
\end{eqnarray*}
with $H_1, \sigma_1$ as in Theorem \ref{main0} and $\sigma_{\rm edge,1}< \infty$ defined in \eqref{sedge}.

\medskip

\noi (ii) Let $Q_{\rm edge,1} <  1  < Q_{\rm edge,2} $. 
Then the convergence in \eqref{partS1} for $\gamma \ne \gamma^\mathfrak{X}_0$ holds
with $\gamma^\mathfrak{X}_0 = \gamma^0_{\rm edge,1},
A_{\la, \gamma} = \lambda^{H(\gamma)} $ and the limit RFs specified in \eqref{V+}, \eqref{V-},
viz.,
\begin{eqnarray*}
\begin{cases} V^\mathfrak{X}_- = \sigma_2 B_{1/2, H_2}, \\
V^\mathfrak{X}_+ = \sigma_{\rm edge,2} B_{0, 1/2},
\end{cases} \qquad
H(\gamma) = \begin{cases}
\gamma H_2 + (1/2), &\gamma < \gamma^0_{\rm edge,1}, \\
\gamma/2, &\gamma > \gamma^0_{\rm edge,1},
\end{cases}
\end{eqnarray*}
with $H_2, \sigma_2$ as in Theorem \ref{main0} and $\sigma_{\rm edge,2}< \infty$ defined in \eqref{sedge}.

\end{thm}

\begin{thm} \label{main2} {\rm [Region $R_{33}$: `two-sided edge effects'.]} Let $Q_{\rm edge,1} \vee Q_{\rm edge,2} < 1$.
Then the convergence in \eqref{partS1} holds with $\gamma^\mathfrak{X}_0 = 1,
A_{\la, \gamma} = \lambda^{H(\gamma)} $ and the limit RFs specified in \eqref{V+}, \eqref{V-},
viz.,
\begin{eqnarray*}
\begin{cases} V^\mathfrak{X}_+ = \sigma_{\rm edge,2} B_{0,1/2}, \\
V^\mathfrak{X}_- = \sigma_{\rm edge,1} B_{1/2,0},  \\
V^\mathfrak{X}_0 = \sigma_{\rm edge,1} B_{1/2,0} + \sigma_{\rm edge,2} B_{0,1/2}
\end{cases} \qquad
H(\gamma) = \begin{cases}
\gamma/2, &\gamma > 1, \\
1/2, &\gamma < 1, \\
1/2, &\gamma = 1,
\end{cases}
\end{eqnarray*}
where $B_{1/2,0} $ and $B_{0,1/2}$ are independent RFs and
$\sigma_{\rm edge, i}, i=1,2$ are given in \eqref{sedge}.

\end{thm}

\begin{remark} \label{rem0} From \eqref{BH0} for any $(x_i, y_i), i=1,2$  we get that
\begin{eqnarray} \label{covBH0}
\E B_{H_1,0}(x_1,y_1) B_{H_1,0}(x_2,y_2)
&=&\frac{1}{2}(x_1^{2H_1}+ x_2^{2H_1} - |x_1-x_2|^{2H_1})\times
\begin{cases}
1,  &y_1 = y_2, \\
\frac{1}{2}, &y_1 \ne y_2,
\end{cases},  \quad 0< H_1 \le  1,  \\
\E B_{0,H_2}(x_1,y_1) B_{0, H_2}(x_2,y_2)
&=&\frac{1}{2}(y_1^{2H_2}+ y_2^{2H_2} - |y_1-y_2|^{2H_2})\times
\begin{cases}
1,  &x_1 = x_2, \\
\frac{1}{2}, &y_1 \ne y_2,
\end{cases}, \quad 0< H_2 \le 1, \nn \\
\E B_{0,0}(x_1,y_1) B_{0,0}(x_2,y_2)
&=&\begin{cases}
1,  &x_1 = x_2, y_1 = y_2,  \\
\frac{1}{2}, &x_1 = x_2, y_1 \ne y_2 \quad \text{or} \quad x_1 \ne x_2, y_1 = y_2, \\
\frac{1}{4}, &x_1 \ne x_2 \quad \text{and} \quad y_1 \ne y_2.
\end{cases} \nn
\end{eqnarray}
Particularly,
\begin{eqnarray*}
\E B_{1/2,0}(x_1,y_1) B_{1/2,0}(x_2,y_2)
&=&(x_1 \wedge x_2)\big(\1(y_1 = y_2) + \frac{1}{2}\1(y_1 \ne y_2)\big),  \\
\E B_{0,1/2}(x_1,y_1) B_{0,1/2}(x_2,y_2)
&=&(y_1 \wedge y_2)\big(\1(x_1 = x_2) + \frac{1}{2}\1(x_1 \ne x_2)\big).
\end{eqnarray*}
The FBS $B_{H_1,H_2}$ with $H_1 \cdot H_2 = 0$ can be defined by Kolmogorov's theorem as Gaussian probability measures on
the space $\R^{\R^2_+}  $ equipped with the $\sigma$-field generated by cylinder sets. They have very irregular (nonseparable) paths
due to the fact that the covariance functions in \eqref{covBH0} are discontinuous  at $(x_1,y_1) = (x_2,y_2)$. Finite dimensional
distributions of these RFs can be easily identified using independent FBM. For instance,
the (joint) distribution of $B_{1/2,0}(x_i,y_i), (x_i,y_i) \in \R^2_+, i=1, \cdots, m, 0< y_1 < \cdots < y_m $ coincides with that of the random vector
\begin{eqnarray}\label{B0m}
\frac{1}{\sqrt{2}}\big( B^{(0)}(x_1) -  B^{(1)}(x_1), B^{(0)}(x_2) -  B^{(2)}(x_2), \cdots, B^{(0)}(x_m) -  B^{(m)}(x_m)\big),
\end{eqnarray}
where $B^{(i)} = \{B^{(i)}(x); x >0\}, i=0,1, \cdots, m$ are {\it independent} standard Brownian motions.  Note
for each $x >0$ fixed, $\{B_{1/2,0}(x,y); y>0 \} \eqfdd \{W(0)- W(y); y >0 \} $  where $W(y) \sim N(0, x/2) \, (\forall \, y\ge 0) $ and
$\E W(y) W(y') = 0 \ (y \ne y', y,y' \ge 0) $.
We note that
the above extension of FBS satisfy the basic (multi) self-similarity property (c.f. \cite{gent2007}): for any
$(H_1,H_2) \in [0,1]^2, \lambda_1 >0, \lambda_2 >0 $
\begin{equation}
B_{H_1,H_2}(\lambda_1^{H_1}x, \lambda_2^{H_2}y) \eqfdd \lambda_1^{H_1} \lambda_2^{H_2} B_{H_1,H_2}(x,y), \quad (x,y) \in \R^2_+.
\end{equation}
\end{remark}

\begin{remark} \label{remE}  The `edge effects' in (i) and (ii) of Theorem \ref{main1} may
be also respectively labelled as `horizontal' and  `vertical'. Note the normalizations $\lambda^{H(\gamma)} = \lambda^{1/2} $ in (i),
$\gamma  < \gamma^0_{\rm edge,2}$ (respectively, $\lambda^{H(\gamma)} = \lambda^{\gamma/2} $ in (ii),
for $\gamma  > \gamma^0_{\rm edge,1}$)
is proportional to
the the square root of the horizontal (respectively, vertical) length of the rectangle $K_{[\lambda x, \lambda^\gamma y]}  $
and $S^\mathfrak{X}_{\lambda, \gamma}(x,y)$ behaves as a sum of weakly dependent r.v.'s  indexed by points on
the horizontal (respectively, vertical) edges of this rectangle. See Remark \ref{edgeremark} below.

\end{remark}

\begin{remark} \label{LRob} In the `isotropic' case $q_1 = q_2 =: \beta $ we have that $Q = 2/\beta \in (0,1)  $ is equivalent to
$\beta \in (2,\infty)$ and $Q_{{\rm edge}, 1} = Q_{{\rm edge}, 2} = 5/(2\beta)$. In this case either Theorem \ref{main0} (for $2< \beta < 5/2 $)
or Theorem \ref{main2} (for $5/2 < \beta < \infty $) applies and our results can be related to (\cite{lah2016}, Theorems 4.1 and 4.3).
In the case $2< \beta < 5/2 $ Theorem \ref{main0} for $\gamma = \gamma^0 = 1$ yields the Gaussian limit
with the variance $\|h_0(x,y;\cdot,\cdot)\|^2 $  which coincides with the variance in (\cite{lah2016}, (4.2)).
In the interval $5/2 < \beta < \infty $ the coincidence of the limiting variances in Theorem \ref{main0}, $\gamma = 1 $ and
(\cite{lah2016}, (4.5)) is less straightforward and can be verified by noting that for any $\delta >0$
\begin{eqnarray}\label{ssigma}
x \sigma^2_{\rm edge,1} + y \sigma^2_{\rm edge,2}
&=&\lim_{\lambda \to \infty} \lambda^{-1} \sum_{(u,v) \in \Z^2: {\rm dist}((u,v), \partial K_{[\lambda x, \lambda y]}) \le \delta \lambda}
\Big(\sum_{(t,s) \in    K_{[\lambda x, \lambda y]}} a(t-u,s-v)\Big)^2,
\end{eqnarray}
where $\partial K_{[\lambda x, \lambda y]} = \{ (t,s)\in  K_{[\lambda x, \lambda y]}: {\rm dist}((t,s), K^c_{[\lambda x, \lambda y]}) = 1\} $
is the boundary of  $K_{[\lambda x, \lambda y]}, K^c_{[\lambda x, \lambda y]} = \Z^2 \setminus K_{[\lambda x, \lambda y]}$ and
$ {\rm dist}((u,v), A) = \inf_{(t,s)\in A} \|(t,s) -  (u,v)\| $ is the distance between $(u,v) \in \Z^2 $ and $A \subset \Z^2 $.
\eqref{ssigma} follows by using the summability properties of $|a(t,s)|$ guaranteeing the convergence
of the series in \eqref{sedge} and the zero-sum condition \eqref{zero}.
\end{remark}

\begin{remark} In Theorem \ref{main1}
the convergence in \eqref{partS1} at $\gamma = \gamma^\mathfrak{X}_0 $  is an open question. We conjecture that
the corresponding limits agree with $B_{1/2,0}$ for $\gamma^\mathfrak{X}_0 = \gamma^0_{\rm edge,2}$ and
$B_{0, 1/2}$ for $\gamma^\mathfrak{X}_0 = \gamma^0_{\rm edge,1}$ however the normalizations might include additional
logarithmic factors making the proofs more difficult. A similar difficulty may accompany the attempts
to extend \eqref{partS1} to values $q_1, q_2 $ belonging to
the boundaries of the partition in Figure~1, particularly to the segments $\widetilde Q_i = 1, Q_{\rm edge,i} = 1, i=1,2 $.
On the other hand, the normalizing exponent $H(\gamma) \equiv H(\gamma, q_1,q_2)$ in Theorems \ref{main0}-\ref{main2} and in
(\cite{pils2017}, Theorems~3.1-3.3, case $k=1$)
extends by continuity to a (jointly) continuous function on the set $\{\gamma >0, 0 < Q < 2 \} $ suggesting
that the above extension of the scaling transition is plausible.

\end{remark}

\begin{remark} Note $\sigma_{\rm edge,1}$ and/or  $\sigma_{\rm edge,2} $
in \eqref{sedge} may vanish provided
\begin{equation}\label{zero1}
\sum_{s \in \Z} a(t,s) = 0  \quad (\forall \ t \in \Z) \qquad \text{and/or} \qquad
\sum_{t \in \Z} a(t,s) = 0  \quad (\forall \ s \in \Z)
\end{equation}
hold. Clearly, each of the conditions in \eqref{zero1} imply \eqref{zero} and the ND property of the corresponding RF, but the converse
implication is not true. Relations \eqref{zero1} may be termed the {\it vertical } (respectively, {\it horizontal}) {\it ND property}.
Under \eqref{zero1} the corresponding limits in Theorems \ref{main1} and \ref{main2} are trivial, while nontrivial limits
may exists under different conditions on $\gamma $ and $q_1, q_2 $, eventually changing Figure~1 and Table~1.

\end{remark}

\section{Proofs of Theorems \ref{main0}-\ref{main2}}

\subsection{Outline of the proof  and preliminaries   
}

Let us explain the main steps of the proof of Theorems \ref{main0}-\ref{main2}. By definition,
$S^\mathfrak{X}_{\la, \gamma}(x,y)$ can be rewritten as
a weighted sum
\begin{equation}\label{SX}
S^\mathfrak{X}_{\la, \gamma}(x,y) = \sum_{(u,v) \in \Z^2} \vep(u,v) G_{\la,\gamma}(u,v),  \qquad
G_{\la,\gamma}(u,v):= \sum_{(t,s) \in K_{[\lambda x, \lambda^\gamma y]}} a(t-u,s-v)
\end{equation}
in i.i.d.r.v.s $\vep(u,v), (u,v) \in \Z^2 $. It happens that different regions of `noise locations' $(u,v)$ contribute to different
limit distributions in our theorems.
\newpage
\begin{center}
\begin{figure}[h]
\begin{center}
\includegraphics[width=15 cm,height=20cm]{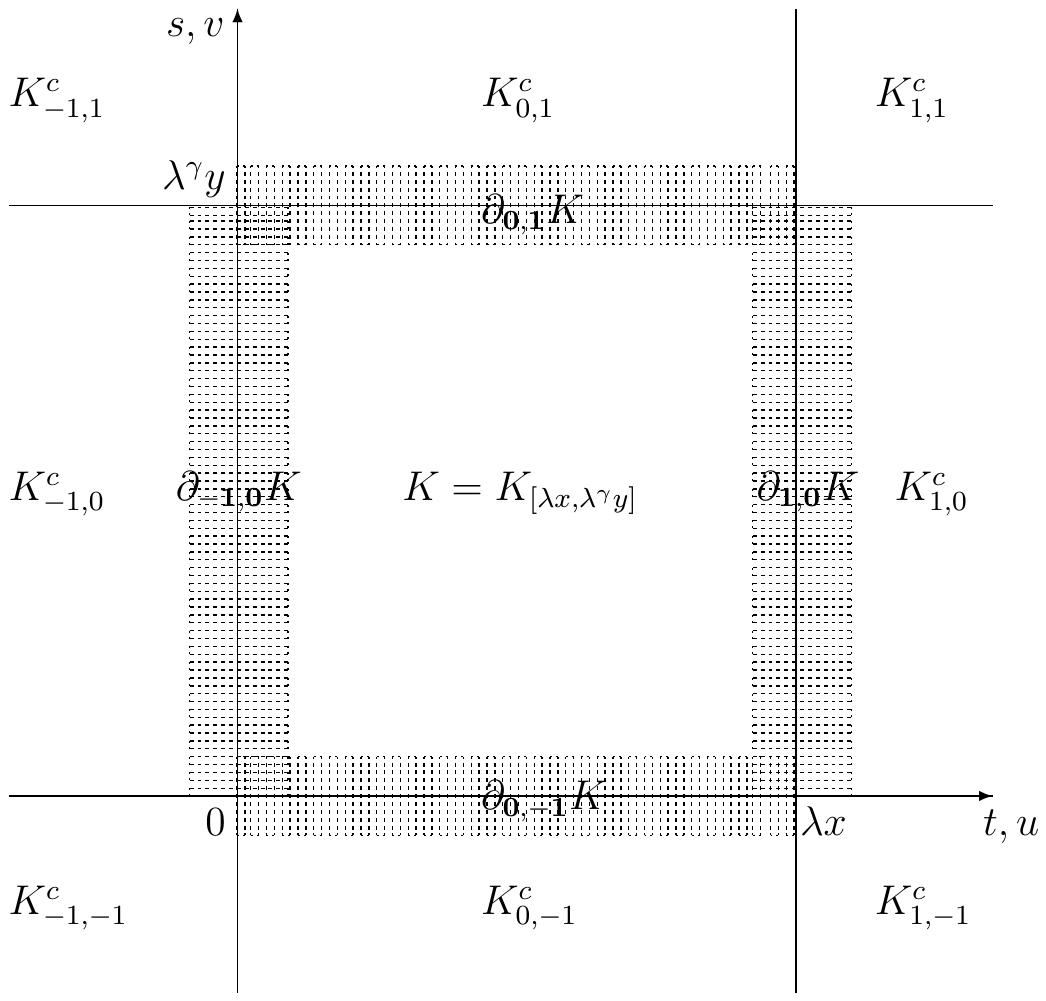}
\vspace{-10.5cm}
\end{center}
\end{figure}
\end{center}

\vskip-2cm

\hskip2cm Figure~2. \parbox[t]{12cm} {\small Partition of $\Z^2 $ by $K = K_{[\lambda x, \lambda^\gamma y]}$ and
8 `outer' sets $
K^c_{i,j}, i,j = 1,0,-1, (i,j) \ne (0,0)$. The shaded regions $\partial_{0,1} K, \partial_{0,-1} K $
(respectively,  $\partial_{1,0} K, \partial_{-1,0} K $)
are related to the `horizontal (respectively, vertical) edge effect'.
}

\bigskip

\noi The basic decomposition
\begin{eqnarray} \label{ZK}
\Z^2 &=&K \cup \big( \bigcup_{(i,j) = -1,0,1, (i,j) \ne (0,0)}  
K^c_{i,j}\big),
\end{eqnarray}
of $\Z^2 $ into 9  sets is shown in Figure~2. Following \eqref{ZK} we decompose \eqref{SX} as
\begin{eqnarray}\label{SK}
&S^\mathfrak{X}_{\la, \gamma}(x,y)\ =  S(K) + \sum_{i,j=-1,0,1, (i,j) \ne (0,0)} S(
K^c_{i,j}),
\end{eqnarray}
where, with zero-sum condition \eqref{zero} in mind,
\begin{eqnarray}\label{SK1}
S(K)&:=&\sum_{(u,v) \in  K} \vep(u,v) G_{\la,\gamma}(u,v) = - \sum_{(u,v) \in  K} \vep(u,v) G^c_{\la,\gamma}(u,v),   \\
S(
K^c_{i,j})&:=&\sum_{(u,v) \in 
K^c_{i,j}} \vep(u,v) G_{\la,\gamma}(u,v),
\qquad G^c_{\la,\gamma}(u,v) := \sum_{(t,s) \not\in K} a(t-u,s-v).  \nn
\end{eqnarray}
As shown  in Lemma  \ref{Krem},
\begin{itemize}

\item terms $S(K), S(K^c_{-1,0}), S(K^c_{1,0}) $ in \eqref{SK}
contribute to the FBS limit $\sigma_1 B_{H_1,1/2} $
for $\gamma > \gamma^\mathfrak{X}_0 $ and
the remaining terms in \eqref{SK} are negligible;

\item terms $S(K), S(K^c_{0,-1}), S(K^c_{0,1}) $ in \eqref{SK}
contribute to the FBS limit $\sigma_2 B_{1/2,H_2}$
 for $\gamma < \gamma^\mathfrak{X}_0 $ and
the remaining terms in \eqref{SK} are negligible;

\item all terms on the r.h.s. of \eqref{SK} are relevant for proving the well-balanced
limit $V^\mathfrak{X}_0$ in Theorem \ref{main0};

\item the main contribution to the unbalanced limit $\sigma_{\rm edge,1} B_{1/2,0}$ in Theorems
\ref{main1}-\ref{main2} comes from innovations `close' to the horizontal edges of $K$ (the shaded regions
$\partial_{0,1} K$ and $\partial_{0,-1} K $ in Figure~2);

\item the main contribution to the unbalanced limit $\sigma_{\rm edge,2} B_{0,1/2}$ in Theorems
\ref{main1}-\ref{main2} comes from innovations `close' to the vertical edges of $K$ (the shaded regions
$\partial_{1,0} K$ and $\partial_{-1,0} K $ in Figure~2).

\end{itemize}

After identification of the main terms, the limit distributions of these terms is obtained using a general criterion
for the convergence of linear forms in i.i.d.r.v.s towards stochastic integral w.r.t. white noise.
Consider a linear form  
\begin{equation} \label{Q}
S(g) := \sum_{(u,v) \in \Z^2} g(u,v) \vep(u,v)
\end{equation}
with real coefficients  $\sum_{(u,v) \in \Z^2} g(u,v)^2 < \infty$.
The following proposition is a version of (\cite{book2012}, Proposition~14.3.2) and
(\cite{sur2019}, Proposition~5.1) and we omit the proof.

\begin{proposition} \label{disc} Let $S(g_\lambda), \lambda >0$ be as in \eqref{Q}. Suppose
$g_\lambda(u,v)$ are such that for a real-valued function
$h \in L^2(\R^2)$ and some 
$m_i = m_i(\lambda) \to \infty, \lambda \to \infty, i=1,2$
the functions
\begin{equation} \label{tildeh}
\tilde g_\lambda(u,v) :=  (m_1 m_2)^{1/2} g_\lambda (\lceil m_1u\rceil, \lceil m_2v\rceil), \quad
(u,v) \in \R^2
\end{equation}
tend to $h$ in $L^2(\R^2)$, viz.,
\begin{equation}\label{L2conv}
\|\tilde g_\lambda - h\|^2 = \int_{\R^2} |\tilde g_\lambda(u,v) - h(u,v)|^2 \d u \d v \to 0, \qquad \lambda \to \infty.
\end{equation}
Then
\begin{equation}\label{Shconv}
S(g_\lambda) \limd I(h) := \int_{\R^2} h(u,v) W(\d u, \d v), \qquad \lambda \to \infty.
\end{equation}

\end{proposition}

\subsection{Two auxiliary lemmas}

The proofs of our results requires evaluation of various multiple series involving
the coefficients $a(t,s)$. Similarly as in \cite{pils2017, sur2019} these series are estimated by corresponding multiple
integrals involving  the function $a_\infty (t,s)$ in \eqref{a0} or
\begin{eqnarray} \label{rho}
&\rho(t,s) \ :=\  (|t|^{q_1} +  |s|^{q_2})^{-1}, \qquad (t,s) \in \R^2. 
\end{eqnarray}
Indeed, by elementary inequalities 
\begin{eqnarray} \label{arho}
&|a(t,s)| \le C \rho(t,s) \quad ((t,s) \in \Z^2), \qquad
|a_\infty (t,s)| \le C \rho(t,s) \quad ((t,s) \in \R^2)
\end{eqnarray}
see (\cite{sur2019}, (2.16)) the function $\rho$ in \eqref{rrho1}-\eqref{aint12}
can be replaced by $|a_\infty|$ while most multiple sums
involving $a(t,s)$ in the proofs below can be evaluated
by corresponding multiple integrals in  Lemma \ref{Lem}.

\begin{lem} \label{Lem} Let $\rho(t,s)$ be defined as in \eqref{rho}, where $q_i>0, i=1,2 $. Then
for any $\delta >0, h >0$
\begin{eqnarray}\label{rrho1}
&&\int_{\R^2} \rho (t,s)^h \1(\rho(t,s) > \delta)
\d t \,\d s < \infty, \qquad \qquad   Q > h, \  \\
\label{rrho2}
&&\int_{\R^2} \rho(t,s)^h \1(\rho(t,s) < \delta)
\d t \, \d s < \infty,  \qquad \qquad Q < h.
\end{eqnarray}
In addition, let $Q < 1 $. Then
\begin{eqnarray}
&&\int_{0}^\infty \d u \Big(\int_1^\infty \d t \int_{0}^\infty \rho(t+u,s)\, \d s  \Big)^2 < \infty,
\qquad \qquad Q_{\rm edge,1} < 1,
\label{aint1} \\
&&\int_{0}^\infty \d u \Big(\int_0^1 \d t \int_{0}^\infty \rho(t+u,s)\, \d s  \Big)^2 < \infty,
\qquad \qquad Q_{\rm edge,1} > 1,
\label{aint1not} \\
&&\int_{0}^\infty \d v \Big(\int_{0}^\infty \d t \int_1^\infty \rho(t,s+v)\, \d s  \Big)^2 < \infty,
\qquad \qquad Q_{\rm edge,2} < 1, \label{aint2} \\
&&\int_{0}^\infty \d v \Big(\int_{0}^\infty \d t \int_0^1 \rho(t,s+v)\, \d s  \Big)^2 < \infty,
\qquad \qquad \ Q_{\rm edge,2} > 1,
\label{aint2not} \\
&&\int_0^\infty \d u  \int_0^1 \d v \Big(\int_0^1 \d t \int_0^1  \rho(t+u,s-v)\, \d s  \Big)^2 < \infty,
\qquad
Q_{\rm edge,1} > 1,
\label{aint3} \\
&&\int_0^1 \d u \int_0^\infty \d v \Big(\int_0^1 \d t \int_0^1 \rho(t-u,s+v)\, \d s  \Big)^2 < \infty,
\qquad Q_{\rm edge,2} > 1,
\label{aint4} \\
&&\int_0^\infty \d u \int_0^\infty \d  v
\Big(\int_0^1 \d t \int_0^1 \rho(t+u,s+v)\, \d s\Big)^2 < \infty, \qquad Q > 2/3,
\label{aint5}
\end{eqnarray}
\begin{eqnarray}
&&\int_0^\infty \d u \int_0^\infty \d  v
\Big(\int_0^1 \d t \int_1^\infty \rho(t+u,s+v)\, \d s\Big)^2 < \infty, \qquad
2/3 < Q < \mbox{$\frac{1}{q_1} + \frac{3}{q_2}$} < 2, \label{aint6} \\
&&\int_0^\infty \d u \int_0^\infty \d  v
\Big(\int_1^\infty \d t \int_0^1 \rho(t+u,s+v) \, \d s\Big)^2 < \infty, \qquad
2/3 < Q < \mbox{$\frac{3}{q_1} + \frac{1}{q_2}$} < 2, \label{aint7} \\
&&\int_{\R^2 \setminus [0,1]^2} \d u \d v \Big(\int_{[0,1]^2} \rho(t-u,s-v)\,\d t \,\d s\Big)^2 < \infty, \qquad
1 < Q_{\rm edge,1}\wedge Q_{\rm edge,2}, \label{aint8} \\
 &&\int_{[0,1]^2} \d u \d v
 \Big(\int_{\R^2 \setminus [0,1]^2} \rho(t-u,s-v)\, \d t \, \d s\Big)^2 < \infty, \qquad
1 < Q_{\rm edge,1}\wedge Q_{\rm edge,2},
\label{aint9}
\end{eqnarray}
Furthermore, as $\mu \to \infty $
\begin{eqnarray}
&&\int_0^\infty \d u \int_0^\infty \d  v
\Big(\int_0^1 \d t \int_1^\mu \rho(t+u,s+v)\, \d s \Big)^2 \nn \\
&&\hskip4cm =  \
O\big(\mu^{3+ \frac{q_2}{q_1}- 2q_2}(1 + \1(q_2 = 3/2) \log \mu)\big), \quad
1 < \mbox{$\frac{1}{2q_1} + \frac{3}{2q_2}$},
\label{aint10} \\
&&\int_0^\infty \d u \int_0^\infty \d  v
\Big(\int_1^\mu \d t \int_0^1 \rho(t+u,s+v)\, \d s\Big)^2 \nn \\
&&\hskip4cm  = \ O\big(\mu^{3 + \frac{q_1}{q_2} - 2q_1}(1 + \1(q_1 = 3/2) \log \mu)\big), \quad
1 < \mbox{$\frac{3}{2q_1} + \frac{1}{2q_2}$};
\label{aint11} \\
&&\int_0^\infty \d  v
\Big(\int_0^1 \d t \int_1^\mu \rho(t,s+v)\, \d s \Big)^2 \   = \
O\big(\mu^{(3-2q_2)\vee 0}\big).
\label{aint12}
\end{eqnarray}

\end{lem}

\smallskip

\noi {\it Proof.} Relations \eqref{rrho1}  and \eqref{rrho2} are proved in (\cite{pils2017}, Proposition~5.1).

\smallskip

\noi Proof of \eqref{aint1}.
Write the l.h.s. of \eqref{aint1} as
$I = \int_0^\infty \d u \big(\int_1^\infty \d t \int_0^\infty
(|t+u|^{q_1} + |s|^{q_2})^{-1} \d s\big)^2 = C  \int_0^\infty \d u \big(\int_1^\infty \d t
(t+u)^{\gamma^0 - q_1} \big)^2 = \int_0^1 \d u \big(\int_1^\infty
t^{\gamma^0 - q_1} \d t\big)^2 + \int_1^\infty \d u (u^{\gamma^0 - q_1 +1})^2 < \infty $
since $2(q_1 -1 - \gamma^0) > 1 $. The proof of \eqref{aint2} is analogous.

\medskip

\noi Proof of \eqref{aint2not}. By integrating the inner integral w.r.t. $t \in (0,\infty)$
the l.h.s. of \eqref{aint2not} can be written as
$I = C  \int_0^\infty \d v \big(\int_0^1
(s+v)^{(q_2/q_1) - q_2}  \, \d s\big)^2 = \int_0^1 \d v (\cdots)^2 + \int_1^\infty \d v (\cdots)^2 =: I_1 + I_2 $ where
$I_1 \le C \int_0^1 v^{2(1 + (q_2/q_1)- q_2)} \d v < \infty $ since $Q_{\rm edge,2} > 1$ and $I_2 \le C \int_1^\infty v^{2((q_2/q_1) - q_2)} \d v < \infty $
since $(1/q_1) + (1/2q_2) < 1 $.
The proof of \eqref{aint1not} is analogous.

\medskip

\noi Proof of \eqref{aint3}. Write  the l.h.s. of \eqref{aint3} as
$I= I_1 + I_2 $ where $I_1 := \int_1^\infty \d u \cdots, I_2 := \int_0^1 \d u  \cdots$.  Then
$I_1 \le C \int_1^\infty \rho(u, 0)^{2} \d u < \infty $ and, by change of variables $t \to u t, s \to u^{q_1/q_2} s $,
$I_2 \le \int_0^1 u^{2 + 2 (q_1/q_2) -2q_1} \d u  < C < \infty $ since
$\int_0^{1/u} \d t\int_0^{1/u^{q_1/q_2}} \rho(t+1, s - (v/u^{q_1/q_2})) \d s
$ $\le C \int_0^\infty \d t \int_{\R} \rho(t+1,s)\, \d s  < \infty $
in view of \eqref{aint2}. The proof of \eqref{aint4} is analogous.

\medskip

\noi   Proof of \eqref{aint5}. Split the l.h.s. of \eqref{aint5}  as
$I= I_1 + I_2$, where  $I_1 := \int_{(0,\infty)^2 \setminus (0,1)^2} \d u \d  v \cdots, I_2 := \int_{(0,1)^2} \d u \d v \cdots$.  Then
$I_1 \le  C\int_{(0,\infty)^2 \setminus (0,1)^2} \rho(u,v)^{2} \d u \d  v \le C $ due to \eqref{rrho2}.
Next by change of variables:  $t \to \rho(u,v)^{-1/q_1} t, s \to \rho(u,v)^{-1/q_2} s $ and using
$\rho(u+ t \rho(u,v)^{-1/q_1}, v + s \rho(u,v)^{-1/q_2}) =   \rho(u,v) \rho(t + u\rho(u,v)^{1/q_1}, s + v \rho(u,v)^{1/q_2}) \ \le
\rho(u,v) (\rho(1,1) \wedge \rho (t,t)), t, s \ge 0 $ we obtain
$$
I_2 \le C \int_{(0,1)^2} \rho(u,v)^{2- 2Q} \d u \d v
\Big(\int_{\R^2} \rho(t,s) \1( \rho(t,s) \le \rho (1,1)) \d t \d s\Big)^2  < \infty
$$
where the first integral converges due to \eqref{rrho1} and
$Q > 2/3$, and the second due to \eqref{rrho2}  and $Q <1$. This proves
\eqref{aint5}.

\medskip

\noi Proof of \eqref{aint6}. Split the l.h.s. of \eqref{aint6}  as
$I= I_1 + I_2$, where  $I_1 := \int_{(0,1) \times \R_+} \d u \d  v \cdots, I_2 := \int_{(1,\infty) \times \R_+} \d u \d v \cdots$.
Then $I_1 \le C \int_0^\infty \d v \big(\int_1^\infty (s+v)^{-q_2} \d s \big)^2
\le C \int_0^\infty \d v \big(1 \wedge v^{-2(q_2 -1)} \big)   < \infty $ since
$2(q_2 -1)  > 1 $ follows from $3/q_2 < 2$.  Next,
$I_2 \le C \int_1^\infty \d u \int_0^\infty \d v \big(\int_1^\infty (u^{q_1} + (s+v)^{q_2})^{-1} \d s\big)^2
\le C J_1 J_2 $, where
$J_1 := \int_1^\infty u^{3(q_1/q_2) - 2q_1} \d u < \infty $ for $\frac{1}{q_1} + \frac{3}{q_2} <2 $ and
\begin{eqnarray*}
J_2&:=&\int_0^\infty \d v  \big(\int_0^\infty (1+ (s+v)^{q_2})^{-1} \d s\big)^2 \\
&\le&C \int_0^1 \d v \big(\int_0^\infty (1+ s^{q_2})^{-1} \d s\big)^2  +
C \int_1^\infty \d v \big(\int_0^v v^{-q_2} \d s  + \int_v^\infty s^{-q_2} \d s \big)^2  \\
&\le&C + C\int_1^\infty v^{2- 2q_2} \d v   < \infty
\end{eqnarray*}
since $q_2 > 3/2$. This proves \eqref{aint6}. The proof of \eqref{aint7} is analogous.

\medskip

\noi Proof of \eqref{aint8}.
Split $\R^2 \setminus [0,1]^2 = \cup_{i=1}^8 D_i $ similarly as in Fig. 2 and then
the l.h.s. of \eqref{aint7} as $I = \sum_{i=1}^8 I_i$ where
$I_i := \int_{D_i} \big(\int_{[0,1]^2} \rho(t-u,s-v) \d t \d s\big)^2  \d u \d v, $  $
D_1 := (-\infty, 0)^2, D_1 := (-\infty, 0)\times (0,1), D_3 := (0,1)\times (-\infty,0)$ etc.
Here, $I_1 < \infty $
according to \eqref{aint5}. 
Similarly, $I_2 =  \int_{(0,\infty)\times (0,1)} \d u  \d v \big(\int_{(0,1]^2} \rho(t+u,s-v) \d t \d s\big)^2 < \infty $
and $I_3 =  \int_{(0,1)\times (0,\infty)} \d u  \d v
\big(\int_{(0,1]^2} \rho(t-u,s+v) \d t \d s\big)^2 < \infty $
follow  from \eqref{aint3} and \eqref{aint4}, respectively.
The remaining relations $I_i < \infty, i=3, \cdots, 8$ follow
from \eqref{aint3}-\eqref{aint5} in a similar way, proving \eqref{aint8}.

\medskip

\noi Proof of \eqref{aint9}. Split $\R^2 \setminus [0,1]^2 = \cup_{i=0}^8 D'_i, $
where $D'_0 := \R^2 \setminus [-1,2]^2, D'_1 := (-1,0)^2, D'_2 := (-1,0)\times (0,1), D'_3 := (-1,0)\times (1,2)$ etc.
Accordingly, for the l.h.s. $J$ of \eqref{aint9} we have $J \le  C \sum_{i=0}^8 J'_i$, where
$J'_i := \int_{[0,1]^2} \d u \d v \big(\int_{D'_i} \rho(t-u,s-v) \d t \d s\big)^2, i=0,1,\cdots, 8$. Note
$\rho(t-u,s-v) \le \delta $ for some $\delta >0$ and
all $(s,v) \in [0,1]^2,
(t,s) \in D'_0 $, implying $J'_0 \le \big(\int_{\R^2} \rho(t,s) \1(\rho(t,s) \le \delta)  \d t \d s\big)^2 < \infty $
according to \eqref{rrho2}.
The proof of $J'_i < \infty, i=1, \cdots, 8$ uses  \eqref{aint3}-\eqref{aint5} and
is analogous to that of $I_i < \infty, i=1, \cdots, 8$ in the proof of  \eqref{aint8}
above. For instance,
$J'_2 =  \int_{[0,1]^2} \d u \d v \big(\int_{(-1,0)\times (0,1)} \rho(t-u,s-v) \d t \d s\big)^2
= \int_{[0,1]^2} \d u \d v \big(\int_{(0,1)^2} \rho(t+u,s-v) \d t \d s\big)^2 < \infty $ according to \eqref{aint3}.
This proves \eqref{aint9}.

\medskip

\noi Proof of \eqref{aint10}. Split the l.h.s. of \eqref{aint10} as
$J_\mu = \int_{\R^2_+ \setminus [0,1]^2} \d u \d v \cdots +
\int_{[0,1]^2} \d u \d v \cdots
=: J'_\mu + J''_\mu. $  Then $J''_\mu \le  \big( \int_1^\infty \rho(0, s) \d s \big)^2 \le C$ since $q_2 > 1$. Next, consider
$J'_\mu = \int_0^1 \d u \int_1^\infty \d v \cdots + \int_1^\infty \d u \int_0^\infty \d v \cdots =: I'_\mu + I''_\mu $. We have
$I'_\mu \le \int_1^\infty \d v \big( \int_1^\mu \rho (0,s+v) \d s\big)^2
= \int_1^\mu \d v \big(\int_1^\mu (s+v)^{-q_2} \d s\big)^2 + \int_\mu^\infty (\mu/v^{q_2})^2 \d v $ and hence
\begin{eqnarray}\label{Iprime}
&I'_\mu  =  O\big(\mu^{(3 - 2q_2)\vee 0}\big) \big(1 + \1(3 = 2q_2) \log \mu\big)
\end{eqnarray}
satisfies the bound in \eqref{aint10}.

Consider
$I''_\mu = \int_1^\infty \d u \int_0^\mu \d v \cdots + \int_1^\infty \d u \int_\mu^\infty \d v \cdots =: L'_\mu + L''_\mu $. Here,
$L''_\mu  \le \int_1^\infty \d u \int_\mu^\infty (\mu \rho(u, v) \big)^2 \d v
= \mu^2 \int_1^\infty \d u $  $ \int_\mu^\infty (u^{q_1} + v^{q_2})^{-2} \d v
= \mu^2 \int_1^\infty u^{(q_1/q_2) - 2q_1} \d u \int_{\mu/u^{q_1/q_2}}^\infty (1+ v^{q_2})^{-2} \d v $
where the last integral is evaluated by splitting it over $ u < \mu^{q_2/q_1} $ and $u > \mu^{q_2/q_1}$ yielding
\begin{eqnarray*}
&L''_\mu  =  O\big(\mu^{3 - 2q_2 + (q_2/q_1) }\big).
\end{eqnarray*}
Finally, $L'_\mu = \int_1^{\mu^{q_2/q_1}} \d u \int_0^\mu \d v \cdots + \int_{\mu^{q_2/q_1}}^\infty \d u \int_0^\mu \d v \cdots
=: L'_{\mu,1} + L'_{\mu,2} $ where
\begin{eqnarray*}
L'_{\mu,1}&\le&\int_1^{\mu^{q_2/q_1}} u^{3(q_1/q_2) -2q_1} \d u
\int_0^{\mu/u^{q_1/q_2}} \d v \big(\int_0^\infty \frac{\d s}{1 + (s+v)^{q_2}} \big)^2   \nn \\
&\le&C\int_1^{\mu^{q_2/q_1}} u^{3(q_1/q_2) -2q_1} \d u
\int_0^{\mu/u^{q_1/q_2}} \frac{\d v}{1 + v^{2(q_2 -1)}}  \ = \  O\big(\mu^{3 - 2q_2 + (q_2/q_1)}(1 + \1(q_2 = 3/2) \log \mu)\big)
\end{eqnarray*}
and
$ L'_{\mu,2}  \le C \mu^3 \int_{\mu^{q_2/q_1}}^\infty u^{-2q_1} \d u
  \ = \  O\big(\mu^{3 - 2q_2 + (q_2/q_1)}\big), $
proving \eqref{aint10}. The proof of \eqref{aint11} is completely analogous, by exchanging $t$ and $s$.

\medskip

\noi Proof of \eqref{aint12}.  Split the l.h.s. of \eqref{aint12} as $I_\mu = \int_0^1 \d v (\cdots)^2 + \int_1^\infty \d v (\cdots)^2
=: I_{\mu,1} + I_{2,\mu}$.  Then $I_{1,\mu} \le \big(\int_1^\infty \rho(0,s) \d s\big)^2 < C $ and
$I_{2,\mu} \le \int_1^\infty \d v \big(\int_1^\mu  (s+v)^{-q_2} \d s\big)^2 $ satisfies the bound in \eqref{aint12} since $q_2 > 1 $.
This proves  \eqref{aint12} and completes the proof of the lemma.
\hfill $\Box$

\medskip

We use the above lemma for evaluation of `remainder terms' in \eqref{SK1}, where
\begin{eqnarray}
K^c_{-1,-1}&:=&\{(u,v) \in \Z^2: u \le 0, v \le 0\}, \quad  K^c_{1,-1} \ :=  \ \{(u,v) \in \Z^2: u > [\lambda x], v \le 0\}, \nn \\
K^c_{0,-1}&:=&\{(u,v) \in \Z^2: 1\le u \le [\lambda x], v \le 0\}, \nn
\end{eqnarray}
the remaining sets analogously defined, viz.,
$K^c_{-1,0} := \{(u,v) \in \Z^2: u \le 0, 1 \le v \le [\lambda^\gamma y]\}, $
$K^c_{1,0} := \{(u,v) \in \Z^2: u > [\lambda x],  1\le v \le [\lambda^\gamma y]\}, $
$K^c_{-1,1} := \{(u,v) \in \Z^2: u \le 0, v > [\lambda^\gamma y]\}, $
$K^c_{0,1} := \{(u,v) \in \Z^2: 1\le u \le [\lambda x], v > [\lambda^\gamma y]\}, $
$K^c_{1,1} := \{(u,v) \in \Z^2: u > [\lambda x], v > [\lambda^\gamma y]\}. $

\begin{lem} \label{Krem} Let  $0 < Q < 1, Q_{\rm edge,i} \ne 1, i=1,2$.  Then
\begin{eqnarray} \label{SKc}
\E \big(S(K^c_{i,j})\big)^2&=&o(\lambda^{2H(\gamma)}), \qquad \lambda \to \infty
\end{eqnarray}
in the following cases:
\begin{eqnarray} \label{corners}
&\gamma > \gamma^\mathfrak{X}_0, 
\quad (i,j) = (-1,-1),\ (1,-1), \ (1,1), \ (-1,1), \ (0,-1), \ (0,1)
\end{eqnarray}
and
\begin{eqnarray} \label{others}
&\gamma < \gamma^\mathfrak{X}_0,
\quad (i,j) = (-1,-1), \ (1,-1), \ (1,1), \ (-1,1), \ (-1,0), \ (1,0).
\end{eqnarray}

\end{lem}

\noi {\it Proof.} In view of the bound in \eqref{arho} and
the reflection symmetry of the function $\rho(t,s)$ w.r.t. the coordinate axes,
it suffices to prove \eqref{corners} for $(i,j) =  (-1,-1) $ and $(i,j) = (-1,0)$. Moreover, it suffices to
consider $x=y=1$.
\smallskip

\noi {\it Case $(i,j) = (-1,-1)$.} Denote $J_{\lambda, \gamma}$ the l.h.s. of
\eqref{SKc} for $x=y=1, i=j= -1$.  Then using \eqref{arho}
\begin{eqnarray}\label{Jbdd}
J_{\lambda, \gamma}
&=&\sum_{u,v \ge 0} \big(\sum_{t=1}^{\lfloor \lambda \rfloor} \sum_{s=1}^{\lfloor \lambda^\gamma \rfloor}
a(t+u,s+v)\Big)^2 \nn \\
&\le&C\int_{\R^2_+} \d u \d v \big(\int_0^\lambda \int_0^{\lambda^\gamma}
\rho(t+u,s+v) \d t  \d s\big)^2 \nn \\
&=&C
\begin{cases}
\lambda^{3 + 3\gamma^0 -2q_1}
\int_{\R^2_+} \d u \d v \big(\int_0^1 \int_0^{\lambda^{\gamma-\gamma^0}}
\rho(t+u,s+v)\, \d t  \d s\big)^2, &\gamma \ge \gamma^0, \\
\lambda^{(3 + 3\gamma^0 -2q_1)(\gamma/\gamma^0)}
\int_{\R^2_+} \d u \d v \big(\int_0^{\lambda^{1- (\gamma/\gamma^0)}} \int_0^{1}
\rho(t+u,s+v)\, \d t  \d s\big)^2,
&\gamma \le \gamma^0.
\end{cases} \label{SK2}
\end{eqnarray}
Then from \eqref{SK2},
\eqref{aint7}, \eqref{aint8}, \eqref{aint10} with $\mu = \lambda^{\gamma - \gamma_0} $  we obtain 
that
\begin{eqnarray} \label{SK3}
J_{\lambda, \gamma}
&\le&C\begin{cases}\lambda^{\phi_+(\gamma)}, &\gamma\ge \gamma^0, \\
\lambda^{\phi_-(\gamma)}, &\gamma \le \gamma^0,
\end{cases}
\qquad \text{where} \nn \\
\phi_+(\gamma)&:=&3 + 3\gamma^0 -2q_1 + (\gamma - \gamma^0)(3  +(1/\gamma^0)  - 2q_2)_+, \nn \\
\phi_-(\gamma)&:=&
\gamma ((3/\gamma^0)  + 3  -2q_2)  + (\gamma^0- \gamma)((3/\gamma^0) +1  - 2q_2)_+. \nn
\end{eqnarray}
Then \eqref{SKc} follows from
\begin{equation}\label{SK4}
2H(\gamma) > \begin{cases}  \phi_+(\gamma), &\gamma > \gamma^0, \\
\phi_-(\gamma), &\gamma < \gamma^0,
\end{cases}
\quad \gamma \ne \gamma^\mathfrak{X}_0,
\end{equation}
with $H(\gamma)$ defined in Theorems \ref{main0}-\ref{main2}. The latter  theorems
determine four of cases of $Q_{{\rm edge},i}, i=1,2 $ in which $\gamma^\mathfrak{X}_0$ takes different values.

\begin{trivlist}

\item [($C$1)] Case $Q_{\rm edge,1} > 1$. Then
$\gamma^\mathfrak{X}_0 \ge  \gamma^0,  2H(\gamma) =
3+ 2\gamma^0 - 2q_1 + \gamma $ (see Theorems \ref{main0} and \ref{main1} (i)) and
$\phi_+(\gamma) = 3 + 3\gamma^0 -2q_1 + (\gamma - \gamma_0)(3  +(1/\gamma^0)  - 2q_2)_+ < 3+ 2\gamma^0 - 2q_1 + \gamma = 2H(\gamma)$
follows by $\gamma > \gamma^0,  Q < 1 $;

\item [($C$2)] Case $Q_{\rm edge,1} <1  < Q_{\rm edge,2}$. Then   $\gamma^\mathfrak{X}_0 = \gamma^0_{\rm edge,1} < \gamma^0,  2H(\gamma) = \gamma $.
For $ \gamma^0_{\rm edge,1} <  \gamma < \gamma^0$  we have that
$((3/\gamma^0) +1  - 2q_2)_+ = 0, \phi_- (\gamma) = \gamma ((3/\gamma^0)  + 3  -2q_2)    $ and
$\phi_- (\gamma) < 2 H(\gamma) $ reduces to $Q_{\rm edge,1} <1$; for $\gamma \ge \gamma^0 $ relation
$\phi_+(\gamma) =3 + 3\gamma^0 -2q_1 + (\gamma - \gamma^0)(3  +(1/\gamma^0)  - 2q_2)_+  < 2H(\gamma) = \gamma $
follows similarly;

\item [($C$3)] Case $Q_{\rm edge,1}\vee Q_{\rm edge,2} < 1$.  Then $\gamma^\mathfrak{X}_0 = 1,  2H(\gamma ) = \gamma  \ (\gamma > \gamma^\mathfrak{X}_0)$
(see Theorem \ref{main2}), $ (3  +(1/\gamma^0)  - 2q_2)_+ = 0,
((3/\gamma^0) +1  - 2q_2)_+ = 0$ and $\phi_+(\gamma)= 3 + 3\gamma^0 -2q_1, $   $
\phi_-(\gamma) =
\gamma ((3/\gamma^0)  + 3  -2q_2) $.  Then $\phi_+ (\gamma) < \gamma \, (\gamma \ge \gamma_0 \vee 1), \,
\phi_- (\gamma ) < \gamma \, (1 < \gamma \le \gamma^0 \vee 1) $ follow
by $Q_{{\rm edge},i} < 1, i=1,2 $.

\end{trivlist}

\smallskip

\noi {\it Case $(i,j) = (0,-1)$.} Write $\widetilde J_{\lambda, \gamma} $ for the l.h.s.
of \eqref{SKc} for the above $i,j$ and $x=y=1$. Then similarly to \eqref{Jbdd}
\begin{eqnarray}
\widetilde J_{\lambda, \gamma}
&\le&C\int_{0}^\lambda  \d u  \int_1^\infty \d v \Big(\int_0^\lambda \int_0^{\lambda^\gamma}
\rho(t-u,s+v)\, \d t  \d s\Big)^2 \nn \\
&\le&C\lambda
\int_0^\infty \d v \Big(\int_0^\lambda \int_1^{\lambda^\gamma}
\rho(t,s+v)\, \d t  \d s\Big)^2 \label{SK5} \\
&=&C
\begin{cases}
\lambda^{3 + 3\gamma^0 -2q_1}
\int_{0}^\infty \d v \big(\int_0^1 \int_0^{\lambda^{\gamma-\gamma^0}}
\rho(t,s+v)\, \d t  \d s\big)^2, &\gamma \ge \gamma^0, \\
\lambda^{1 + (2 + 3\gamma^0 -2q_1)(\gamma/\gamma^0)}
\int_{0}^\infty \d v \big(\int_0^{\lambda^{1- (\gamma/\gamma^0)}} \int_0^{1}
\rho(t,s+v)\, \d t  \d s\big)^2,
&\gamma \le \gamma^0.
\end{cases} \nn \\
&\le&C\begin{cases}\lambda^{\tilde \phi_+(\gamma)}, &\gamma \ge \gamma^0, \\
\lambda^{\tilde \phi_-(\gamma)}, &\gamma \le \gamma^0,
\end{cases}
\label{SK6}
\end{eqnarray}
where
\begin{eqnarray} \label{SK7}
\tilde \phi_+(\gamma)&:=&3 + 3\gamma^0 -2q_1  + (\gamma-\gamma^0)(3-2q_2)_+,   \\
\tilde \phi_-(\gamma)&:=&\begin{cases}1, & Q_{\rm edge,2} < 1, \\
1 + (2 + 3\gamma^0 -2q_1)(\gamma/\gamma^0), &Q_{\rm edge,2} > 1 > (1/q_1)  + (1/2q_2)
\end{cases}   \nn
\end{eqnarray}
follow from \eqref{SK5}, \eqref{SK6} and \eqref{aint2}, \eqref{aint2not}.
Then \eqref{SKc} for $(i,j) = (0,-1)$ follows from
\begin{equation}\label{SK8}
2H(\gamma) > \begin{cases} \widetilde \phi_+(\gamma), &\gamma > \gamma^0, \\
\widetilde \phi_-(\gamma), &\gamma < \gamma^0,
\end{cases}
\quad \gamma \ne \gamma^\mathfrak{X}_0,
\end{equation}
with $H(\gamma)$ defined in Theorems \ref{main0}-\ref{main2}.

\begin{trivlist} 

\item [($\widetilde C$1)] Case $Q_{\rm edge,1} > 1$. Then
$\gamma^\mathfrak{X}_0 \ge  \gamma^0,  2H(\gamma) =
3+ 2\gamma^0 - 2q_1 + \gamma $  and
$\tilde \phi_+(\gamma) = 3 + 3\gamma^0 -2q_1 + (\gamma - \gamma_0)(3 - 2q_2)_+ < 2H(\gamma)$ as in ($C$1);

\item [($\widetilde C$2)] Case $Q_{\rm edge,1} <1  < Q_{\rm edge,2}$. Then   $\gamma^\mathfrak{X}_0 = \gamma^0_{\rm edge,1} < \gamma^0,  2H(\gamma) = \gamma $.
For $ \gamma^0_{\rm edge,1} <  \gamma < \gamma^0$  we have that
$\widetilde \phi_- (\gamma) =  1 + (2 + 3\gamma^0 -2q_1)(\gamma/\gamma^0) <   \gamma
= 2H(\gamma) $ is equivalent to  $\gamma > \gamma^0_{\rm edge,1} $ (condition  $1 > (1/q_1)  + (1/2q_2)$ in \eqref{SK7} is satisfied
by $Q_{\rm edge,1} <1$); for $\gamma \ge \gamma^0$ relation $\widetilde \phi_+ (\gamma) < \gamma $ follows
by $Q_{\rm edge,1} <1$;

\item [($\widetilde C$3)] Case $Q_{\rm edge,1}\vee Q_{\rm edge,2} < 1$.  Then $\gamma^\mathfrak{X}_0 = 1,  2H(\gamma ) = \gamma  \ (\gamma > \gamma^\mathfrak{X}_0)$.
Let $\gamma > \gamma^0 \vee 1 $ then $\widetilde \phi_+ (\gamma) < \gamma $ as in ($\widetilde C$2). Next, let $1 < \gamma \le \gamma^0$ then
$\widetilde \phi_- (\gamma) = 1 < \gamma = 2H(\gamma)$.

\end{trivlist}

\noi The above discussion proves \eqref{SKc} for \eqref{corners}.  The proof of \eqref{SKc} for
\eqref{others} is analogous and omitted.  Lemma \ref{Krem} is proved. \hfill $\Box$

\medskip

\subsection{Proof of Theorem \ref{main0}}

The finite-dimensional convergence in
\eqref{partS1}, or $ A^{-1}_{\la, \gamma} S^\mathfrak{X}_{\la,\gamma}(x,y) \limfdd V^\mathfrak{X}_\gamma (x,y)$,
towards the corresponding limit RF $V^\mathfrak{X}_\gamma$
is equivalent to the one-dimensional convergence of arbitrary linear combinations, viz.,
$A^{-1}_{\la, \gamma} \sum_{\ell=1}^p \theta_\ell S^\mathfrak{X}_{\la,\gamma}(x_\ell,y_\ell) $  $ \limd \sum_{\ell=1}^p V^\mathfrak{X}_\gamma (x_\ell, y_\ell)$
for any  $\theta_\ell \in \R, (x_\ell, y_\ell) \in \R^2_+, \ell =1, \cdots, p, $  $ p \ge 1. $
Towards this aim we use Proposition  \ref{disc} and write $A^{-1}_{\la, \gamma} \sum_{\ell=1}^p \theta_\ell S^\mathfrak{X}_{\la,\gamma}(x_\ell,y_\ell)
=  Q(g_{\lambda,\gamma}) $ as a linear form in \eqref{Q} with coefficients $g_{\lambda,\gamma}(u,v)= \sum_{\ell =1}^p \theta_\ell g_{\lambda,\gamma}(x_\ell,y_\ell; u,v)$ where
\begin{equation} \label{glambda}
g_{\lambda,\gamma}(x,y; u,v) \ := \  \lambda^{-H(\gamma)}
\sum_{(t,s)\in  K_{[\lambda x,\lambda^{\gamma} y]}}  a(t-u,s-v), \quad
(u,v) \in \Z^{2}.
\end{equation}
Since in all cases of Theorem \ref{main0} the limit RF $V^\mathfrak{X}_\gamma (x,y) = \int_{\R^2} h(x,y; u,v) W (\d u, \d v)$ writes
as a white noise integral w.r.t. some $h(x,y; \cdot, \cdot) \in L^2(\R^2)$, by   Proposition  \ref{disc} it
it suffices to prove the $L^2$-convergence
\begin{equation} \label{ghconv}
\|\tilde g_{\lambda,\gamma} - h\|^2  \to 0
\end{equation}
with $h(u,v) :=   \sum_{\ell=1}^p \theta_\ell h_\gamma (x_\ell,y_\ell; u,v)$
and suitably chosen $m_i = m_i(\lambda), i=1,2 $, see \eqref{tildeh}, \eqref{L2conv}. Clearly,
the proof of \eqref{ghconv} can be reduced to the case $p=\theta_1 = 1 $ and $(x_1,y_1) = (x,y) \in \R^2_+ $ fixed.
Accordingly, the proof of \eqref{partS1} is reduced to one-dimensional convergence
$ A^{-1}_{\la, \gamma} S^\mathfrak{X}_{\la,\gamma}(x,y) \limd V^\mathfrak{X}_\gamma (x,y)$
at a given point $(x,y) \in \R^2_+$, or the verification of \eqref{ghconv} for
$g_{\lambda, \gamma}(\cdot, \cdot) \equiv g_{\lambda,\gamma}(x,y; \cdot, \cdot)$ in \eqref{glambda} and
$h(\cdot,\cdot) \equiv h_j (x,y; \cdot, \cdot), j=0,1,2$
in \eqref{h00LRD}
in respective cases $\gamma = \gamma^0, \gamma > \gamma^0, \gamma < \gamma^0$.
For this, it is convenient
to rewrite $g_{\lambda,\gamma}$ and $\tilde g_{\lambda,\gamma}$ as integrals of step functions
and change the variables $(t,s) \to (\tilde m_1 t, \tilde m_2 s): \R^2 \to \R^2$, with $\tilde m_i, i=1,2 $ defined below.
Then using \eqref{zero}
\begin{eqnarray}\label{tildeg}
\tilde g_{\lambda,\gamma}(x,y; u,v)
&=&\kappa_{\lambda,\gamma}
\begin{cases}
\int_{\widetilde K_{\lambda, \gamma}(x,y)}
a(\lceil \tilde m_1 t \rceil - \lceil m_1 u\rceil, \lceil \tilde m_2 s \rceil - \lceil m_2 v \rceil) \d t \d s,
&(u,v) \not\in  \widetilde K_{\lambda, \gamma}(x,y), \\
-\int_{\widetilde K^c_{\lambda, \gamma}(x,y)}
a(\lceil \tilde m_1 t \rceil - \lceil m_1 u\rceil, \lceil \tilde m_2 s \rceil - \lceil m_2 v \rceil) \d t \d s,
&(u,v) \in  \widetilde K_{\lambda, \gamma}(x,y),
\end{cases}
\end{eqnarray}
where
\begin{eqnarray}\label{wtig1}
&\widetilde K_{\lambda,\gamma}(x,y):= \big(0, \lfloor \lambda x \rfloor/\tilde m_1 \big] \times
\big(0, \lfloor \lambda^\gamma y \rfloor/\tilde m_2 \big], \quad
\widetilde K^c_{\lambda,\gamma}(x,y) := \R^2\setminus  \tilde K_{\lambda,\gamma}(x,y),  \\
&\kappa_{\lambda, \gamma}:= (m_1 m_2)^{1/2} \tilde m_1  \tilde m_2 \, \lambda^{-H(\gamma)}.  \nn
\end{eqnarray}
The normalising factor $\kappa_{\lambda, \gamma}$ will be used to prove the limit of the integrand on the r.h.s. of \eqref{tildeg}.

Consider \eqref{ghconv} for $\gamma = \gamma^0$. Let $m_1 = \tilde m_1 := \lambda,
m_2 =  \tilde m_2 := \lambda^{\gamma^0}. $ Then
$\widetilde K_{\lambda,\gamma^0}(x,y) = \big(0, \lfloor \lambda x \rfloor/\lambda \big] \times
 \big(0, \lfloor \lambda^{\gamma^0} y \rfloor/\lambda^{\gamma^0} \big],  \kappa_{\lambda, \gamma^0}=  \lambda^{q_1} $
and
\begin{eqnarray}\label{tildegN}
\tilde g_{\lambda,\gamma^0}(x,y; u,v)
&=&
\begin{cases}
\int_{\widetilde K_{\lambda, \gamma^0}(x,y)}
a_\lambda (t,s,u,v) \d t \d s,
&(u,v) \not\in  \widetilde K_{\lambda, \gamma^0}(x,y), \\
-\int_{\widetilde K^c_{\lambda, \gamma^0}(x,y)}
a_\lambda(t,s,u,v) \d t \d s,
&(u,v) \in  \widetilde K_{\lambda, \gamma^0}(x,y),
\end{cases}
\end{eqnarray}
where
\begin{eqnarray}
a_\lambda (t,s,u,v)&:=&\lambda^{q_1} a((\lceil \lambda t \rceil - \lceil \lambda u\rceil,
\lceil \lambda^{\gamma^0} s \rceil - \lceil \lambda^{\gamma^0} v \rceil)  \to a_\infty(t-u,s-v)
\end{eqnarray}
for any  $(t,s) \ne (s,v)$. Moreover,
\begin{eqnarray} \label{ala}
|a_\lambda(t,s,u,v)| \le C\rho(t-u,s-v), \qquad (t,s) \in \R^2, \ (s,v) \in \R^2
\end{eqnarray}
see (\cite{pils2017}, (7.28)), with $C >0$ independent of $\lambda \ge 1 $ and $\rho $ as in \eqref{rho}. Whence, the point-wise convergence
$\tilde g_{\lambda,\gamma_0}(x,y; u,v)  \to h_{0} (x,y;  u,v), (u,v) \in \R^2, (u,v) \not \in  \partial K(x,y)$ (= the boundary of
$K(x,y) := (0,x]\times (0,y]$)
easily follows. Moreover, since $  \widetilde K_{\lambda, \gamma^0}(x,y) \subset K(x,y)$ we see from
\eqref{ala} that $\tilde g_{\lambda,\gamma^0}(x,y; u,v) \le \bar h(u,v) + \bar h_\lambda(u,v)$ where
$$
\bar h_0(x,y;u,v) := C\begin{cases}
\int_{K(x,y)} \rho(t-u,s-v) \d t \d s,  &(u,v) \not\in K(x,y), \\
\int_{K^c(x,y)} \rho(t-u,s-v) \d t \d s,  &(u,v) \in K(x,y),
\end{cases}
$$
and $\bar h_\lambda(u,v):= \1((u,v) \in K(x,y)\setminus \widetilde K_{\lambda, \gamma^0}(x,y)) \int_{\widetilde K^c_{\lambda, \gamma^0}(x,y)} \rho(t-u,s-v) \d t \d s.  $   Then $\| \bar h \| < \infty, \|\bar h_\lambda \| \to 0 $ follow from Lemma \ref{Lem},
\eqref{aint8}-\eqref{aint9}. The above facts together with the dominated convergence theorem
prove \eqref{ghconv} for $j=0$ or $\gamma = \gamma^0$.
\smallskip

Next, consider \eqref{ghconv} for $j=1$ or $\gamma > \gamma^0$.
Let
$m_1 = \tilde m_1 := \lambda,
m_2 := \lambda^{\gamma}, \tilde m_2 := \lambda^{\gamma^0} $.  Using Lemma \ref{Krem}, in \eqref{ghconv}
we can replace $\tilde g_{\lambda, \gamma}(x,y;\cdot,\cdot)$ of \eqref{tildeg}  by its restriction
on $\R \times (0, \lfloor \lambda^\gamma y \rfloor/\lambda^\gamma]$, or  the function
\begin{eqnarray}\label{tildeg1}
h_{\lambda,\gamma}(x,y; u,v)
&:=&\1(0< v \le \lfloor \lambda^\gamma y \rfloor/\lambda^\gamma)
\begin{cases}
\int_{(0,\lfloor \lambda x\rfloor/\lambda]\times \R} a_{\lambda, \gamma}(t,s,u,v) \d t \d s,
&u \in  (0, \lfloor \lambda x \rfloor/\lambda], \\
-\int_{\R^2 \setminus (0,\lfloor \lambda x\rfloor/\lambda]  \times \R}
a_{\lambda,\gamma}(t,s,u,v) \d t \d s,
&u \not\in  (0, \lfloor \lambda x \rfloor/\lambda],
\end{cases}
\end{eqnarray}
where
\begin{eqnarray*}
a_{\lambda, \gamma}(t,s,u,v)&:=&\kappa_{\lambda,\gamma}
a(\lceil \lambda t \rceil - \lceil \lambda u\rceil, \lceil \lambda^{\gamma^0} s \rceil)
\1( -\lceil \lambda^\gamma v \rceil/\lambda^{\gamma^0} < s \le \lceil \lfloor \lambda^\gamma y \rfloor
- \lambda^\gamma v \rceil/\lambda^{\gamma^0})
\end{eqnarray*}
tends to $a_\infty(t-u,s) $ point-wise for any $t \ne u, s \in \R, 0 < v <y $. Hence
as $ -\lceil \lambda^\gamma v \rceil/\lambda^{\gamma^0} \to - \infty, \lceil \lfloor \lambda^\gamma y \rfloor
- \lambda^\gamma v \rceil/\lambda^{\gamma^0} \to  \infty \ (0 < v < y) $ we see that $h_{\lambda,\gamma}(x,y; u,v) \to h_1(x,y;u,v) $
in \eqref{h00LRD} point-wise for any $(u,v) \in \R^2, u \not\in \{0, x\}, v \not\in \{0,y\}$.
Then the $L^2(\R^2)$-convergence in
 \eqref{ghconv} for $j=1$ follows similarly as in the case $j=0$ above.
The proof of \eqref{ghconv} for $j=2$ or $\gamma < \gamma^0$
is similar using $m_1 = \lambda, \tilde m_1 =   \lambda^{\gamma/\gamma^0},  m_2 = \tilde m_2 = \lambda^\gamma $ and
Lemma \ref{Krem}.  Theorem \ref{main0} is proved. \hfill $\Box$

\subsection{Proof  of Theorem \ref{main1}}

(i) 
The
convergence in \eqref{partS1} for $\gamma  > \gamma^0_{\rm edge,2}$
follows from the the proof of
Theorem \ref{main0} or the convergence in \eqref{ghconv} with
$\tilde g_{\lambda, \gamma}(x,y;\cdot,\cdot)$ replaced by
$h_{\lambda,\gamma}(x,y; u,v)$ of
\eqref{tildeg1}.

Let us prove \eqref{partS1} for $\gamma  < \gamma^0_{\rm edge, 2}$. Using  Lemma \ref{Krem}
we can there replace $S^\mathfrak{X}_{\lambda,\gamma}(x,y)$ by
\begin{eqnarray*}
S^{0}_{\lambda,\gamma}(x,y)&:=&\sum_{1\le u \le  \lfloor\lambda x\rfloor, v \in \Z}
\vep(u,v) \sum_{(t,s) \in K_{[\lambda x,\lambda^{\gamma} y]}}  a(t-u,s-v).
\end{eqnarray*}
Split  $S^0_{\lambda,\gamma}(x,y)= R^0_{\lambda,\gamma}(x,y) + R^+_{\lambda,\gamma}(x,y)+ R^-_{\lambda,\gamma}(x,y)$ where
\begin{eqnarray*}
R^{-}_{\lambda,\gamma}(x,y)&:=&\sum_{1 \le u \le  \lfloor\lambda x\rfloor,\, v \le 0} \vep(u,v)
\sum_{(t,s) \in K_{[\lambda x,\lambda^{\gamma} y]}}  a(t-u,s-v), \\
R^{0}_{\lambda,\gamma}(x,y)&:=&\sum_{1\le u \le  \lfloor\lambda x\rfloor, \, 1 \le v \le  \lfloor\lambda^\gamma y\rfloor }
\vep(u,v) \sum_{(t,s) \in K_{[\lambda x,\lambda^{\gamma} y]}}  a(t-u,s-v),  \\
R^{+}_{\lambda,\gamma}(x,y)&:=&\sum_{1 \le u \le  \lfloor\lambda x],\, v >  \lfloor\lambda^\gamma y\rfloor}
\vep(u,v) \sum_{(t,s) \in K_{[\lambda x,\lambda^{\gamma} y]}}  a(t-u,s-v).
\end{eqnarray*}
where $R^{0}_{\lambda,\gamma}(x,y)$ is further rearranged as
\begin{eqnarray*}
R^{0}_{\lambda,\gamma}(x,y)&:=&\sum_{(u,v) \in K_{[\lambda x,\lambda^{\gamma} y]}}
\vep(u,v) \sum_{(t,s) \in K_{[\lambda x,\lambda^{\gamma} y]}}  a(t-u,s-v)  \\
&=&-\sum_{(u,v) \in K_{[\lambda x,\lambda^{\gamma} y]}}
\vep(u,v) \sum_{(t,s) \not\in K_{[\lambda x,\lambda^{\gamma} y]}}  a(t-u,s-v) \\
&=&M^0_{\lambda,\gamma}(x,y) + M^+_{\lambda,\gamma}(x,y)+ M^-_{\lambda,\gamma}(x,y)
\end{eqnarray*}
where
\begin{eqnarray*}
M^{-}_{\lambda,\gamma}(x,y)&:=&-\sum_{1\le u \le \lfloor \lambda x\rfloor, \, 1 \le v \le  \lfloor\lambda^\gamma y\rfloor} \vep(u,v)
\sum_{t \le 0, 1 \le s \le  \lfloor\lambda^\gamma y\rfloor} a(t-u,s-v), \\
M^{0}_{\lambda,\gamma}(x,y)&:=&-\sum_{1\le u \le  \lfloor\lambda x\rfloor,\, 1 \le v \le  \lfloor\lambda^\gamma y\rfloor}
\vep(u,v) \sum_{1 \le t \le  \lfloor\lambda x\rfloor,\, s \not\in [1,  \lfloor\lambda^\gamma y\rfloor]}  a(t-u,s-v),  \\
M^{+}_{\lambda,\gamma}(x,y)&:=&-\sum_{1 \le u \le  \lfloor\lambda x\rfloor,\, 1\le v  \le  \lfloor\lambda^\gamma y\rfloor}
\vep(u,v) \sum_{t >  \lfloor\lambda x\rfloor, \, 1\le v \le  \lfloor\lambda^{\gamma} y\rfloor}  a(t-u,s-v).
\end{eqnarray*}
Let us prove that  $M^{\pm}_{\lambda,\gamma}(x,y)$ are negligible, viz.,
\begin{eqnarray}\label{M0bdd}
\E (M^\pm_{\lambda,\gamma}(x,y))^2
&=&\sum_{1 \le u \le  \lfloor\lambda x\rfloor,\, 1 \le v \le  \lfloor\lambda^\gamma y\rfloor}
\big(\sum_{t \le 0, 1 \le s \le  \lfloor\lambda^{\gamma} y\rfloor }  a(t-u,s-v)\big)^2
\ = \ o(\lambda), \quad \gamma < \gamma^0_{\rm edge,2}.
\end{eqnarray}
\eqref{M0bdd} follows by integral approximation as in the proof of Lemma \ref{Krem}. Indeed, using \eqref{arho} the l.h.s. of
\eqref{M0bdd} can be evaluated as
\begin{eqnarray*}\label{J}
&&C\int_1^\lambda \d u \int_0^{\lambda^\gamma} \d v \big(\int_0^\infty \d t \int_0^{\lambda^\gamma}
\rho(t+u,s-v) \d s\big)^2 \
\le\ C \lambda^{\gamma} \int_1^\lambda \d u \big(\int_0^\infty \d t \int_0^{\infty}
\rho(t+u,s) \d s\big)^2 \\
&&=\ C \lambda^{\gamma} \int_1^\lambda \d u \big(\int_0^\infty (t+u)^{(q_1/q_2) - q_1} \d t\big)^2
\ = \ C \lambda^{\gamma} \int_1^\lambda  u^{2 (1 + (q_1/q_2)  -q_1)} \d u \\
&&= \ C \lambda^{\gamma + 3 + 2(q_1/q_2)-2q_1 } \
= \ o(\lambda), \quad \text{for} \quad  \gamma < \gamma^0_{\rm edge,2}.
\end{eqnarray*}
We finally arrive at the main term:
\begin{eqnarray}
M_{\lambda,\gamma}(x,y)&:=&R^{-}_{\lambda,\gamma}(x,y) + R^{+}_{\lambda,\gamma}(x,y) + M^{0}_{\lambda,\gamma}(x,y) \nn \\
&=&\sum_{1 \le u \le \lfloor \lambda x\rfloor, \, v \not\in [1, \lfloor \lambda^\gamma y \rfloor]} \vep(u,v)
\sum_{1 \le t \le \lfloor \lambda x\rfloor, \, s \in [1, \lfloor \lambda^\gamma y\rfloor ]}  a(t-u,s-v) \nn \\
&-&\sum_{1\le u \le \lfloor \lambda x\rfloor, \, 1 \le v \le \lfloor \lambda^\gamma y \rfloor}
\vep(u,v) \sum_{1 \le t \le \lfloor \lambda x\rfloor, \, s \not\in [1, \lfloor \lambda^\gamma y\rfloor]}  a(t-u,s-v)\nn \\
&=&\sum_{1 \le u, t \le \lfloor \lambda x\rfloor } \Big\{\sum_{s \in [1, \lfloor \lambda^\gamma y\rfloor],\, v \not\in
[1,\lfloor \lambda^\gamma y\rfloor ]}
-  \sum_{s \not\in [1, \lfloor \lambda^\gamma y\rfloor],\, v \in [1, \lfloor \lambda^\gamma y\rfloor]} \Big\}
\vep(u,v) a(t-u,s-v). \label{Main}
\end{eqnarray}
Let us show the  convergence
\begin{equation}\label{Mainconv}
\lambda^{-1/2} M_{\lambda,\gamma}(x,y)
\ \limfdd\
\sigma_{\rm edge,1} B(x), \quad \gamma < \gamma^0_{\rm edge,2}.
\end{equation}
Note $|s-v| \ge 1 $ in \eqref{Main} and $\min \{|s-v|; s,v \in \Z\} $  is achieved
at the lower and upper rows $ \{1 \le t \le \lfloor \lambda x\rfloor, s=1 \}  $ and
$\{1 \le t \le \lfloor \lambda x\rfloor, s = \lfloor \lambda^\gamma y\rfloor\} $ of the rectangle  $K_{[\lambda x,\lambda^{\gamma} y]} $.
For a large but fixed $L \ge 1 $ decompose the l.h.s. of \eqref{Mainconv} as
$M_{\lambda,\gamma}(x,y) = \sum_{i=0}^4 V^{(i)}_{L,\lambda,\gamma}(x,y)$, where
\begin{eqnarray*}
V^{(1)}_{L,\lambda,\gamma}(x,y)
&:=&\sum_{1 \le u, t \le \lfloor \lambda x\rfloor }\ \sum_{-L <  s \le 0 < v \le L} \vep(u,v)   a(t-u,s-v), \\
V^{(2)}_{L,\lambda,\gamma}(x,y)
&:=&-\sum_{1 \le u, t \le \lfloor \lambda x\rfloor}\ \sum_{-L <  v \le 0 < s \le L} \vep(u,v) a(t-u,s-v), \\
V^{(3)}_{L,\lambda,\gamma}(x,y)
&:=&-\sum_{1 \le u, t \le \lfloor \lambda x\rfloor}\  \sum_{-L <  s - \lfloor \lambda^\gamma y\rfloor  \le 0 < v - \lfloor \lambda^\gamma y\rfloor  \le L} \vep(u,v)
a(t-u,s-v), \\
V^{(4)}_{L,\lambda,\gamma}(x,y)
&:=&\sum_{1 \le u, t \le \lfloor \lambda x\rfloor}\  \sum_{-L <  v - \lfloor\lambda^\gamma y\rfloor \le 0 < s - \lfloor \lambda^\gamma y\rfloor  \le L } \vep(u,v)
a(t-u,s-v),
\end{eqnarray*}
$V^{(0)}_{L,\lambda,\gamma}(x,y) := M_{\lambda,\gamma}(x,y) - \sum_{i=1}^4 V^{(i)}_{L,\lambda,\gamma}(x,y)$.
Note $V^{(i)}_{L,\lambda,\gamma}(x,y), i=1,2,3,4 $ are independent and have similar structure,
each consisting of a finite number $L$ of  (independent) weighted sums of noise variables $\vep(u,v), 1 \le t \le \lfloor \lambda x\rfloor $
in a neighborhood of the lower ($v=0$) or upper ($v = \lceil \lambda^\gamma y\rceil $) boundaries of the
rectangle   $K_{[\lambda x,\lambda^{\gamma} y]} $; moreover, the distribution of these sums does not depend on $y$. We claim  that
\begin{eqnarray}\label{V1}
&&\lambda^{-1/2} V^{(i)}_{L,\lambda,\gamma}(x,y)
\ \limfdd\
\sigma^{(i)}_{L} B(x), \quad i=1,2,3,4, \  1\le L < \infty,  \\
&&\lim_{L\to \infty} \limsup_{\lambda \to \infty} \lambda^{-1} \E \big(V^{(0)}_{L,\lambda,\gamma}(x,y))^2 \ = \ 0, \label{V2} \\
&&\lim_{L \to \infty} \sum_{i=1}^4 (\sigma^{(i)}_{L})^2 \ = \ \sigma^2_{\rm edge,1}. \label{V3}
\end{eqnarray}
The asymptotic variances in \eqref{V1} are given by
\begin{eqnarray}\label{S1}
(\sigma^{(i)}_{L})^2&:=&\begin{cases}
\sum_{v=0}^{L-1} \big(\sum_{t \in \Z, 1\le s \le L} a(t,s+v)\big)^2, & i=2,4, \\
\sum_{v=1}^L \big(\sum_{t \in \Z, 0\le s < L} a(t,-s-v)\big)^2, & i=1,3.
\end{cases}
\end{eqnarray}
We omit the proof of \eqref{V1}-\eqref{V3} which reduces to a standard application of Lindeberg's theorem and
the dominated convergence theorem and relies on the fact that the series for $\sigma^2_{\rm edge, 1}$
in \eqref{sedge} absolutely converges:
\begin{eqnarray}\label{sedge2}
&\sum_{v \ge 0} \big(\sum_{t\in \Z, s \ge 1} |a(t,s+v)| \big)^2 +
\sum_{v \le -1} \big(\sum_{t\in \Z, s \le 0} |a(t,s+v)| \big)^2 \ < \ \infty.
\end{eqnarray}
The convergence in \eqref{sedge2} follows from \eqref{aint2} and the assumption $Q_{\rm edge,2} < 1 $. By  Slutsky's theorem,
\eqref{V1}-\eqref{V3} imply \eqref{Mainconv}, thereby completing the proof of Theorem \ref{main1} (i).
The proof of  Theorem \ref{main1} (ii) is analogous and is omitted. \hfill $\Box$

\subsection{Proof  of Theorem \ref{main2}}

For $\gamma < \gamma^\mathfrak{X}_0 = 1$ the convergence in \eqref{Mainconv} and the approximation
$S^\mathfrak{X}_{\lambda,\gamma}(x,y)= M_{\lambda,\gamma}(x,y) + o_p(\lambda^{1/2})$ follow as in the proof
of Theorem \ref{main1} (i); particularly,  \eqref{M0bdd} holds for $\gamma < 1 $ in view
of \eqref{aint1} and $Q_{\rm edge,1} < 1 $. For $\gamma > 1 $ the proof of the theorem is analogous using the approximation
$S^\mathfrak{X}_{\lambda,\gamma}(x,y)= N_{\lambda,\gamma}(x,y) + o_p(\lambda^{\gamma/2})$, where
\begin{eqnarray*}
N_{\lambda,\gamma}(x,y)
&:=&\sum_{1 \le v, s \le \lfloor \lambda^\gamma y\rfloor} \Big\{\sum_{t \in [1, \lfloor \lambda x\rfloor], u \not\in [1, \lfloor \lambda x\rfloor]}
-  \sum_{t \not\in [1, \lfloor \lambda x\rfloor], u \in [1, \lfloor \lambda x\rfloor]} \Big\}
\vep(u,v) a(t-u,s-v) 
\end{eqnarray*}
is the analog of $M_{\lambda,\gamma}(x,y)$ in \eqref{Main} and satisfies
\begin{equation}\label{Nainconv}
\lambda^{-\gamma/2} N_{\lambda,\gamma}(x,y)
\ \limfdd\
\sigma_{\rm edge,2} B_2(y) 
\end{equation}
where $B_2 $ is a standard Brownian motion. Finally, for $ \gamma = 1$ we have the  approximation
$S^\mathfrak{X}_{\lambda,\gamma}(x,y)= M_{\lambda,\gamma}(x,y) +
N_{\lambda,\gamma}(x,y) + o_p(\lambda^{1/2})$ and the joint convergence
in \eqref{Mainconv} and \eqref{Nainconv} towards independent BM $B = B_1$ and $B_2$, thus proving
Theorem \ref{main2}.  \hfill  $\Box$

\begin{remark}\label{edgeremark}  It follows from the proofs of Theorems  \ref{main1}--\ref{main2} that the FBS 
limits $B_{1/2,0} $ and $B_{0,1/2}$ in Table~1 (associated with the `edge effects' in Remark \ref{remE})
essentially
reduce to the approximation of
$S^\mathfrak{X}_{\lambda,\gamma}(x,y)$ by sums of horizontal and/or vertical  `edge processes'   defined below.  Particularly,
\begin{eqnarray}\label{SU}
&\E \Big(S^\mathfrak{X}_{\lambda,\gamma}(x,y) - \sum_{1 \le t \le \lfloor \lambda x\rfloor} (U_{0}(t) - U_{\lfloor \lambda^\gamma  y\rfloor}(t))\Big)^2
= o(\lambda^{2H(\gamma)}), \quad \text{for} \quad \gamma < \gamma^\mathfrak{X}_0,
\end{eqnarray}
where,  for each $s \in \Z $
the `horizontal  process' $U_{s} = \{U_{s}(t); t \in \Z \} $ `located' on the line $\{(t,s) \in \Z^2: t \in \Z \} \subset \Z^2$ is defined by
\begin{eqnarray}
&U_{s}(t):=
\sum_{(u,v,w) \in \Z^3: {\rm sgn}(w- s - (1/2)) \neq  {\rm sgn}(v- s - (1/2)) } a(t-u,w-v) {\rm sgn}(v-w)   \vep(u, v),  \quad t \in \Z.
\end{eqnarray}
The process $U_s$ is stationary and SRD which follows from
\begin{eqnarray*}
\sum_{t \in \Z} |{\rm Cov}(U_s(0), U_s(t))|
&\le&\sum_{v \in \Z}\big( \sum_{ (u,w) \in \Z^2} |a(u,w-v)| \1({\rm sgn}(w-(1/2)) \ne {\rm sgn}(v-(1/2))) \big)^2
< \infty,
\end{eqnarray*}
see  \eqref{sedge2}. 
Moreover, $\sum_{t \in \Z} {\rm Cov}(U_s(0), U_s(t)) = (1/2) \sigma^2_{\rm edge,1} $ and the
 `edge processes'  $ U_{0}$ and $ U_{\lfloor \lambda^\gamma y\rfloor}$ in \eqref{SU} `located' on the lower and upper edges
 of the rectangle $K_{[\lambda x, \lambda^\gamma y]}$
 become uncorrelated as $\lambda \to \infty $. 
The above facts may explain the asymptotic variance  $\sigma^2_{\rm edge,1} $ in \eqref{sedge} by approximation in \eqref{SU}. The
convergence $\lambda^{-1/2} \sum_{1\le t \le \lfloor \lambda x\rfloor } U_s(t) \limfdd (1/\sqrt{2}) \sigma_{\rm edge,1} B(x) $ follows
similarly as in \eqref{V1}.
Analogous approximation of $S^\mathfrak{X}_{\lambda,\gamma}(x,y)$ by vertical `edge processes' holds
in Theorems \ref{main1}-\ref{main2} for
$\gamma > \gamma^\mathfrak{X}_0$.
\end{remark}

\begin{remark}\label{remB0} As noted above,  FBS $B_{1/2,0} $ and $B_{0,1/2}$
may arise under anisotropic scaling of very simple finitely-dependent RF with ND. Consider a linear RF $\mathfrak{X}$ in \eqref{Xlin} with
coefficients
\begin{equation}\label{a}
a(t, s) := \1(t=0, s = 0) - \1(t=0, s=1), \qquad (t,s) \in \Z^2.
\end{equation}
Then $S^{\mathfrak{X}}_{\lambda,\gamma}(x,y) =  \sum_{1 \le t \le \lfloor \lambda x\rfloor} \vep(t,1) - \sum_{1 \le t \le \lfloor \lambda x\rfloor} \vep(t,
\lfloor \lambda^\gamma y\rfloor +1)  $ is the difference of two independent sums of i.i.d.r.v.s, each of which tends
to a standard BM under normalization  $\lambda^{1/2}$. It easily follows that for each  $y>0$ fixed,
$\lambda^{-1/2} S^{\mathfrak{X}}_{\lambda,\gamma}(x,y) \limfdd  B^{(0)}(x) - B^{(1)}(x)$ where $B^{(0)}, B^{(1)} $ are independent BM. In a similar way,
for any $(x_i,y_i) \in \R^2_+, i=1, \cdots, m$ as in \eqref{B0m} we have the joint convergence in distribution of
the random vector $(2\lambda)^{-1/2} \big(S^{\mathfrak{X}}_{\lambda,\gamma}(x_i,y_i), i= 1, \cdots, m\big) $ towards the Gaussian
vector in \eqref{B0m} implying the convergence $(2\lambda)^{-1/2}S^{\mathfrak{X}}_{\lambda,\gamma}(x,y) \limfdd B_{1/2,0}(x,y)$
for any $\gamma >0$.  Note that \eqref{a} satisfy the vertical ND property in  \eqref{zero1}.

\end{remark}

\section{Fractionally integrated negative dependent RFs}

\noi {\it 1. Isotropic fractionally integrated RF.}
Introduce the (discrete) Laplacian
$\Delta Y(t,s) := (1/4) \sum_{|u|+|v|=1} $  $(Y(t+u,s+v) - Y(t,s))$
and a lattice isotropic linear 
RF
\begin{equation}\label{Xfrac1}
\mathfrak{X}(t,s)\  = \
\sum_{(u,v)\in \Z^2} a(u,v) \vep(t-u,s-v), \qquad (t,s) \in \Z^2
\end{equation}
where $\{\vep(t,s), (t,s) \in \Z^2 \}$ are standard i.i.d. r.v.s,
\begin{eqnarray}
a(u,v)&:=&\sum_{j=0}^\infty \psi_j(-d) p_j(u,v), \qquad
\psi_j(d) := \Gamma(j-d)/\Gamma(j+1) \Gamma(- d), \qquad |d| < 1/2
\end{eqnarray}
and $p_j (u,v) $ 
are $j$-step transition probabilities  of a symmetric  nearest-neighbor random walk $\{W_j; j=0,1,\cdots \}$
on $\Z^2$ with equal 1-step probabilities $\P(W_1 = (u,v)|W_0 = (0,0)) =
1/4, |u| + |v| = 1$. Note $\sum_{j=0}^\infty |\psi_j(d)| < \infty,  \sum_{j=0}^\infty \psi_j(d) = 0 \ (-1/2 < d < 0), \
\sum_{j=0}^\infty \psi_j(d)^2 < \infty \ (0< d < 1/2) $ implying
\begin{eqnarray} \label{asum}
&&\sum_{(t,s) \in \Z^2} |a(t,s)| < \infty, \qquad \sum_{(t,s) \in \Z^2} a(t,s) = 0 \qquad (-1/2 < d < 0), \\
&&\sum_{(t,s) \in \Z^2} a(t,s)^2 < \infty  \qquad (0 < d < 1/2), \nn
\end{eqnarray}
see also \cite{kou2016, pils2017}. As shown in the latter papers,
for $ 0< d < 1/2 $ the linear RF $\mathfrak{X}$ in  \eqref{Xfrac1} is the unique stationary solution of the fractional equation
\begin{equation}\label{Xfrac12}
(- \Delta)^d \mathfrak{X}(t,s) \ = \  \vep(t,s),
\end{equation}
where the operator on the l.h.s. is defined as
\begin{eqnarray}\label{Xfrac11}
(- \Delta)^{d} Y(t,s)
&=&\sum_{j=0}^\infty \psi_j(d) (1+\Delta)^j Y(t,s)  \
=\  \sum_{(u,v) \in \Z^2} b(u,v) Y(t-u,s-v),  \\
b(u,v)&:=&\sum_{j=0}^\infty \psi_j(d) p_j(u,v), \quad (u,v) \in \Z^2. \nn
\end{eqnarray}
Moreover, for $ 0< d < 1/2 $ the  RF $\mathfrak{X}$ in  \eqref{Xfrac1} is LRD
with moving-average coefficients satisfying \eqref{acoefL} with $q_1 = q_2 = 2(1-d) \in (1,2)$ and a constant angular function $L_0$.
By stationary solution of \eqref{Xfrac12} we mean a covariance stationary RF $\mathfrak{Y} = \{\mathfrak{Y}(t,s); (t,s) \in \Z^2 \} $, with finite second moment
$\E \mathfrak{Y}(t,s)^2 < \infty $
such that the series $(- \Delta)^{d} \mathfrak{Y}(t,s)
=\sum_{(u,v) \in \Z^2} b(u,v) \mathfrak{Y}(t-u,s-v) $ converges in mean square  and $ (- \Delta)^{d} \mathfrak{Y}(t,s) =  \vep(t,s) $ holds,
for each $(t,s) \in \Z^2 $.
The following proposition extends the above-mentioned result to negative $d \in (-1/2,0)$.

\begin{prop} \label{propB} Let $-1/2 < d < 0$. Then:

\smallskip

\noi (i) RF $\mathfrak{X}$ in \eqref{Xfrac1} is a stationary solution of equation \eqref{Xfrac11}. Moreover,
this solution is unique among the class of all linear RFs $ \mathfrak{Y}(t,s) = \mu + \sum_{(u,v) \in \Z^2} c(u,v) \vep (t-u,s-v), (t,s) \in \Z^2 $
with $\sum_{(u,v) \in \Z^2} |c(u,v)| < \infty $ and $\mu = \E \mathfrak{Y}(t,s) = 0$.

\smallskip

\noi (ii) $\mathfrak{X}$ in \eqref{Xfrac1}
is a ND RF and
\begin{equation} \label{ado}
a(t,s) = (A + o(1))(t^2 + s^2)^{-(1- d)}, \quad t^2 + s^2 \to \infty,
\end{equation}
where $A := \pi^{-1} \Gamma (1- d)/\Gamma (d) <  0.$ Particularly, $a(t,s)$ satisfy \eqref{acoefL} with $q_1 = q_2 = 2(1-d) \in (2,3), Q = 1/(1-d) \in (2/3,1)$
and a constant angular function $L_0(z) = A$.

\end{prop}

\noi {\it Proof.} We use the spectral representation $\vep (t,s) = \int_{\Pi^2} \e^{\i (tx+sy)} Z(\d x, \d y), (t,s) \in \Z^2, $ where
$\Pi := [-\pi,\pi], Z(\d x, \d y)$ is a (random) complex-valued spectral measure, $\overline{Z(\d x, \d y)} = Z(-\d x, -\d y),
\E Z(\d x, \d y) = 0,   \E |Z(\d x, \d y)|^2 = (2\pi)^{-2} \d x \d y.  $ Then $\mathfrak{X}(t,s)$ in \eqref{Xfrac1} can be written as
\begin{eqnarray}\label{Xspec}
&\mathfrak{X}(t,s) = \int_{\Pi^2}  \e^{\i (tx+sy)} \widehat a(x,y) Z(\d x, \d y),
\end{eqnarray}
where   $\widehat a(x,y) $ is the Fourier transform:
\begin{eqnarray}
&\widehat a(x,y) := \sum_{(u,v)\in \Z^2}    \e^{-\i (ux+vy)} a(u,v) = (1 - \widehat p_1(x,y))^{-d},
\end{eqnarray}
$\widehat p_1(x,y) := (1/4) \sum_{|u|+|v|=1}    \e^{-\i (ux+vy)} = (\cos x + \cos y)/2$,
see (\cite{kou2016}, (5.5)). Applying the operator  $(-\Delta)^d $
\eqref{Xfrac11} to $X $ in \eqref{Xspec} we obtain
$(-\Delta)^d \mathfrak{X}(t,s) = \int_{\Pi^2} \e^{\i (tx+sy)} \big(\sum_{(u,v) \in \Z^2}  \e^{-\i (ux+vy)} b(u,v) \big) \widehat a(x,y) Z(\d x, \d y)
= \int_{\Pi^2} \e^{\i (tx+sy)} \widehat b(x,y)  \widehat a(x,y) Z(\d x, \d y) =  \vep(t,s) $ since
$\widehat b(x,y) = 1/\widehat a(x,y), $ the series $ \sum_{(u,v) \in \Z^2}  \e^{-\i (ux+vy)} b(u,v) = \widehat b(x,y)$ converges
in $L^2(\Pi^2)$ and $\widehat a(x,y)$ is bounded on $\Pi^2 $. To show the uniqueness part, let $ \mathfrak{Y}(t,s) = $  $
\sum_{(u,v) \in \Z^2} c(u,v) \vep (t-u,s-v)$ be another solution having spectral representation
$\mathfrak{Y}(t,s) = \int_{\Pi^2}  \e^{\i (tx+sy)} \widehat c(x,y)$  $ Z(\d x, \d y)$ where $\widehat c(x,y)$ is a bounded function on $\Pi^2 $. Then \
$(-\Delta)^d (\mathfrak{X}(t,s) - \mathfrak{Y}(t,s)) = \int_{\Pi^2} \e^{\i (tx+sy)} \widehat b(x,y)$  $ (\widehat a(x,y)) - \widehat c(x,y)) Z(\d x, \d y) = 0$ or
$\int_{\Pi^2} |\widehat b(x,y)|^2 |\widehat c(x,y)- \widehat a(x,y)|^2 \d x\,\d y  = 0$ implying $ \widehat c(x,y)=  \widehat a(x,y) $ a.e. on
$\Pi^2 $ since $|\widehat b(x,y)|^2 = |1 - \widehat p_1(x,y)|^{2d} >0 $ a.e. on $\Pi^2 $. Therefore, $c(t,s) = a(t,s) $ and
$\mathfrak{Y}(t,s) = \mathfrak{X}(t,s), (t,s) \in \Z^2$.

\smallskip

\noi (ii) The ND property of $\mathfrak{X}$ \eqref{Xfrac1} follows from the zero-sum condition in \eqref{asum}.
The proof of \eqref{ado} using the Moivre-Laplace theorem carries over from (\cite{kou2016}, proof of Proposition~5.1) to
the  case $-1/2 < d <0$ virtually  without any change. \hfill $\Box$

\begin{remark} It follows from \eqref{Xspec} that
RF $\mathfrak{X}$ has an explicit
spectral density $f(x,y) = (2\pi)^{-2} |\widehat a(x,y)|^2 =
 2^{-2d}|(1-\cos x) +
(1-\cos y)|^{-2 d},  (x,y) \in \Pi^2$ which vanishes as ${\rm const} \, (x^2 + y^2)^{-2d} \to  0$ as
$x^2 +  y^2 \to 0$.

\end{remark}

\noi{\it 2. Anisotropic ND fractionally integrated RF.} Following \cite{pils2017}
consider the `discrete heat operator'
$\Delta_{1,2} \mathfrak{X}(t,s)$  $ := \mathfrak{X}(t,s) -  \theta \mathfrak{X}(t-1, s) - \frac{1-\theta}{2} (\mathfrak{X}(t-1,s+1) + \mathfrak{X}(t-1,s-1)), 0< \theta < 1$
and a corresponding `fractional power' defined as
\begin{equation}\label{Xfrac21}
\Delta_{1,2}^{d} \mathfrak{X}(t,s) \ := \  \sum_{(u,v)\in \Z^2} b(u,v) \mathfrak{X}(t-u,s-v),
\end{equation}
where
\begin{equation}\label{buv}
b(u,v) = \psi_u(d) q_\theta(u,v)\1(u \ge 0),  \qquad -3/4 < d < 0,
\end{equation}
where  $q_\theta(u,v) $
are $u$-step transition probabilities  of a random walk $\{W_u; u=0,1,\cdots \}$
on $\Z$ with 1-step probabilities $\P(W_1 = v|W_0 = 0) =
\theta $ if $v = 0$, $=  (1- \theta)/2$ if $v= \pm 1$. As shown in \cite{pils2017}, $\sum_{(u,v) \in \Z^2} b(u,v)^2 < \infty $
and
\begin{eqnarray}\label{L0}
b(u,v)&=&\frac{1}{(|u|^2 + |v|^{2q_2/q_1})^{q_1/2}}\Big(L_0\big(\frac{u}{(|u|^2 + |v|^{2q_2/q_1})^{1/2}}; d, \theta\big) + o(1)\Big), \quad
|u|+|v| \to \infty,
\end{eqnarray}
where
\begin{eqnarray} \label{L1}
&&q_1 = (3/2) + d,   \qquad q_2 =  2q_1  = 2((3/2) + d),  \\
&&L_0(z; d, \theta) \ = \ \begin{cases}
\frac{z^{-d- (3/2)}}{\Gamma (-d) \sqrt{2 \pi (1-\theta)}}  \exp \big\{ - \frac{\sqrt{(1/z)^2 -1}}{2(1-\theta)}\big\},
&0 < z \le 1, \\
0, & -1 \le z \le 0
\end{cases} \nn
\end{eqnarray}
and $L_0(z; d, \theta), z \in [-1,1] $ is a bounded continuous function.
Similarly to \eqref{Xfrac1},
the inverse operator can be defined by
$\Delta_{1,2}^{-d} \mathfrak{X}(t,s) \ := \  \sum_{(u,v)\in \Z_+ \times \Z} a(u,v) \mathfrak{X}(t-u,s-v), $ where
\begin{equation}\label{auv}
a(u,v) = \psi_u(-d) q_\theta(u,v)\1(u \ge 0),  \qquad -3/4 < d < 0,
\end{equation}
are obtained from $b(u,v)$ in \eqref{buv} with $d$ replaced by $-d$. By direct inspection of the proof of
Proposition~4.1 in \cite{pils2017} we find that $a(t,s)$ also satisfy \eqref{L0}-\eqref{L1} with $d $ replaced by $-d$, for $-3/4 < d < 0$, in other words,
the linear RF $\mathfrak{X}$ in \eqref{Xfrac1} with moving-average coefficients in \eqref{auv} satisfies the assumptions in \eqref{acoefL} with
$q_1 = (3/2)-d, q_2 =  2q_1  = 2((3/2) - d) $ and
$$
Q = \frac{2}{3-2d} + \frac{1}{3- 2d} \in (\frac{2}{3}, 1), \qquad -\frac{3}{4} < d < 0.
$$
Moreover, the above $\mathfrak{X}$  can be regarded as a stationary solution (in the sense as explained before Proposition \ref{propB})
of the fractional equation:
\begin{equation}\label{Xfrac22}
\Delta_{1,2}^{d} \mathfrak{X}(t,s) \ = \  \vep(t,s),  \qquad (t,s) \in \Z^2.
\end{equation}
We arrive at the following proposition whose proof is similar to that of Proposition \ref{propB} and
we omit the details.

\begin{prop} \label{propC} Let $-3/4 < d < 0$. Then:

\smallskip

\noi (i) RF $\mathfrak{X}$ in \eqref{Xfrac1} with coefficients $a(u,v)$ in \eqref{auv}
is a stationary solution of equation \eqref{Xfrac22}. Moreover,
this solution is unique among the class of all linear RFs $ \mathfrak{Y}(t,s) = \mu + \sum_{(u,v) \in \Z^2} c(u,v) \vep (t-u,s-v), (t,s) \in \Z^2 $
with $\sum_{(u,v) \in \Z^2} |c(u,v)| < \infty $ and $\mu = \E \mathfrak{Y}(t,s) = 0$.

\smallskip

\noi (ii) The above  $\mathfrak{X}$
is a ND RF and $a(u,v) $ satisfy \eqref{acoefL} with $q_1 = (3/2)-d, q_2 =  2q_1  = 2((3/2) - d) $ and a continuous
angular function $L_0(z) = L_0(z; \theta, -d), z \in [-1,1] $ as defined in \eqref{L1}.

\end{prop}

\noi{\it 3. Separately fractionally integrated RF.} Several works discussed RFs where the operators of
fractional integration are
applied separately w.r.t. each of the two coordinates on the plane. The moving-average coefficients
of such RFs have a product form which is very different from \eqref{acoefL} and apparently does not allow for
scaling transition.
See (\cite{pils2017}, Remark~4.1), \cite{dam2017}, (\cite{sur2019}, Remark~4.1). Gaussian LRD RFs
whose spectral density has a product form appear in \cite{leo2013, ps2015}.
Particularly,
\cite{bois2005, guo2009} discussed fractionally integrated RFs  satisfying the equation
\begin{equation}\label{Xfrac3}
\nabla_1^{d_1} \nabla_2^{d_2}  \mathfrak{X}(t,s) \ = \  \vep(t,s), \qquad (t,s) \in \Z^2,
\end{equation}
where $\nabla_1 \mathfrak{X}(t,s) :=  \mathfrak{X}(t,s) - \mathfrak{X}(t-1,s), \nabla_2 \mathfrak{X}(t,s) := \mathfrak{X}(t,s) - \mathfrak{X}(t,s-1) $ are difference operators
and $0<d_1, d_2 < 1/2 $ are parameters. Stationary solution of \eqref{Xfrac3} is a moving-average RF in
$\mathfrak{X}(t,s) = \sum_{(u,v) \in \Z_+^2} a(u,v)\vep(t-u,v-s)$ with coefficients
$a(u,v) := \psi_u(-d_1) \psi_v (-d_2).$  Equation \eqref{Xfrac3} and its solution  can be extended to any $ d_i \in (-1/2, 1/2), i=1,2 $.
Since the covariance function of this  RF factorizes into the product of two covariances
of ARFIMA processes with respective parameters $d_1, d_2 $, it follows that
the RF $\mathfrak{X}$ in  \eqref{Xfrac3} is LRD if $d_1 \vee d_2 >0$ and ND if $d_1\vee d_2 \le 0, d_1 \wedge d_2 < 0$.
Moreover, following (\cite{pils2017}, Theorem 3.1)  and Theorem \ref{main0} of the present paper,
one can show that for any $\gamma >0$ and any $|d_i| < 1/2, i=1,2 $
 the (normalized) partial sums process of RF $\mathfrak{X}$ in  \eqref{Xfrac3} tends to a FBS depending on $d_1, d_2$ only, viz.,
$\lambda^{-H_1 - \gamma H_2} S^\mathfrak{X}_{\lambda,\gamma}(x,y) \limfdd c(d_1) c(d_2) B_{H_1,H_2}(x,y)  $, where
$H_i = d_i + 1/2 $ and $c(d_i) >0$ are some constants. As a consequence, the separately fractionally integrated RF in \eqref{Xfrac3}
does not exhibit scaling transition in both (LRD and ND) parameter regions.

\cite{espe2015} discussed separately fractionally integrated {\it Gegenbauer} RF satisfying
\begin{equation}\label{Xfrac4}
\nabla_{u_1}^{d_1}  \nabla_{u_2}^{d_2} \mathfrak{X}(t,s)
 \ = \  \vep(t,s), \qquad (t,s) \in \Z^2,
\end{equation}
where $\nabla_{u_i}^{d_i} :=  (1- 2u_i(1- \nabla_i) + (1-\nabla_i)^2)^{d_i},
|u_i| < 1,  d_i \in (0,1/2), i  =1,2 $ are parameters.
\cite{espe2015} proved that $\mathfrak{X}$ in  \eqref{Xfrac4} is a LRD RF which exhibits seasonal or cyclical dependence, due to the fact
that its spectral density is positive and continuous
at the origin and is singular on two lines parallel to the coordinate axes.
We conjecture that all scaling limits of $\mathfrak{X}$ in  \eqref{Xfrac4} agree with the Brownian Sheet $B_{1/2,1/2} $ and
the above RF does not exhibit scaling transition. See \cite{ould2003} for related results in the time series setting.

\begin{remark} \cite{leo2011} discuss a wide class of
fractionally integrated Gaussian RFs on $\R^\nu, \nu \ge 1 $, defined
via Fourier transform, together with numerous references.
Scaling properties of these RFs seem to be less
investigated.

\end{remark}

\section{Concluding remarks}

The present paper and \cite{pils2017} provide a nearly complete description,  summarised in Table~1,
of anisotropic scaling limits on inflated rectangles
of linear LRD, SRD and ND
RFs on $\Z^2 $ with coefficients decaying at generally different rate in the horizontal and vertical directions.
It is proved that under ND these limits may exhibit horizontal and/or vertical edge effects, meaning that the main contribution
to the limit comes from values of the RF near the horizontal and/or vertical edges of the rectangle.
We expect that these results can be generalized to linear RFs in higher dimensions $\nu \ge 3 $, as well
as to linear RFs with infinite variance, although the description of the anisotropic scaling limits
in \cite{sur2019} for LRD  RFs
and $\nu=3 $  is more complicated. 
We also expect that our results  remain valid, in some sense,
for more general sampling regions similar to \cite{lah2016} and  `anisotropically inflated' by scaling factors $\lambda $ and $\lambda^\gamma $
in the horizontal  and vertical directions.

The scaling results in Table~1 strongly rely on scaling properties of moving-average coefficients $a(t,s) $ in \eqref{acoefL}.
See (\cite{sur2019}, Remark~4.1).  In the Gaussian case ($\vep(0,0) \sim N(0,1) $) the RF $\mathfrak{X}$  \eqref{Xlin} is completely
described by its spectral density (which exists for $ Q< 2 $ due to square summability of $a(t,s)$, see \cite{pils2017}).
For Gaussian LRD RFs on $\Z^2 $  similar scaling results were obtained in \cite{ps2015}, assuming that the spectral density
has an asymptotic form at the origin  which resembles that of the moving-average coefficients at infinity. A rigorous
statement relating \eqref{acoefL} to the behavior of its Fourier transform at the origin (Abelian/Tauberian theorem \cite{leo2013b}) is an interesting open problem; see  (\cite{sur2019}, Remark~2.2).

A direction on the plane (a line passing through the origin) along which the moving-average coefficients $a(t,s) $
decay at the smallest rate may be
called the {\it dependence axis} of linear RF $\mathfrak{X}$ in \eqref{Xlin}. Consider the coefficients $a(t,s) $ as in
\eqref{acoefL} where $L_0 \equiv 1 $, say. Then $a(t,0) \sim |t|^{-q_1}, a(0,s) \sim |s|^{-q_2} $ and
the dependence axis of the corresponding RF $\mathfrak{X}$ agrees with the horizontal axis if $q_1 < q_2$ and with
the vertical axis if $q_1 > q_2$.
Since the scaling in \eqref{partS1} is parallel to the coordinate axes,
we may say that for RF $\mathfrak{X}$ in \eqref{Xlin}-\eqref{acoefL} {\it the scaling is congruous with the dependence axis of $\mathfrak{X}$}
and the results in  Table~1 refer to this situation. This raises the question what happens  when
the dependence axis does {\it not} agree with
the coordinate axes (the case of {\it incongruous scaling}), particularly, the question
about the scaling transition and the scaling transition point $\gamma^\mathfrak{X}_0$  under incongruous scaling.
A linear RF with `oblique' dependence axis can be obtained by a `rotation' of $a(t,s) $ in \eqref{acoefL}
or a more general linear transformation of the plane.
We conjecture that the fact of `oblique' dependence axis or
incongruous scaling dramatically changes the scaling transition
so that the scaling transition point is always 1: $\gamma^{\mathfrak{X}}_0 =1 $  independently of $q_1, q_2 $,
and the unbalanced limits are generally different from
the corresponding limits in the congruous scaling case (in the LRD region $1< Q < 2$ this conjecture was recently
confirmed in \cite{pils2020}).

Finally, we expect that many results in Table~1 can be carried over to RFs indexed by continuous argument $(t,s) \in \R^2 $.
Particularly, this remark applies to stationary Gaussian RFs which are completely described by spectral  measure on $\R^2 $.
A popular class of nongaussian RFs with continuous argument is {\it shot-noise},
see  \cite{alb1994}, \cite{gir1992} and the references therein. However, the scaling diagram of  LRD shot-noise RF may
be very different from Table~1. See \cite{pils2016}. These and further open problems are discussed in more detail in the recent review paper
\cite{sur2019b}.

\section*{Acknowledgments}

The author is grateful to two anonymous referees for useful comments.

\hskip2cm

\bigskip


\end{document}